\documentclass{article}
\usepackage{a4,amssymb,amsmath,color,graphicx}
\usepackage{trfsigns}
\usepackage{babel}
\usepackage{siunitx}

\def\eins{\mbox{1\hskip-0.24em l}}
\def\T{^{\sf T}}
\def\mT{^{\sf -T}}

\newcommand{\C}{ {\mathbb C} }

\newcommand{\R}{ {\mathbb R} }
\newcommand{\MM}{{\mathbb M}}

\newcommand{\CC}{{\cal C}}

\newcommand{\PP}{{\cal P}}
\newcommand{\QQ}{{\cal Q}}
\newcommand{\Sc}{{\cal S}}
\newcommand{\diag}{\,\mbox{diag}}
\newcommand{\bldiag}{\,\mbox{blockdiag}}

\newcommand{\cc}{{\bf c}}

\newcommand{\folgt}{\ \Rightarrow\ }

\newcommand{\qed}{\qquad\mbox{$\square$}}
\newcommand{\mynegspace}{\hspace{-0.12em}}
\newcommand{\ltnorm}{\rvert\mynegspace\rvert\mynegspace\rvert}
\newcommand{\rtnorm}{\rvert\mynegspace\rvert\mynegspace\rvert}
\DeclareMathOperator*{\argmin}{arg\,min}

\newtheorem{remark}{Remark}[section]
\newtheorem{theorem}{Theorem}[section]
\newtheorem{lemma}{Lemma}[section]

\begin{document}
\title{Implicit Third-Order Peer Triplets\\
with Variable Stepsizes for Gradient-Based Solutions\\
in Large-Scale ODE-Constrained Optimal Control}
\author{Jens Lang \\
{\small \it Technical University Darmstadt,
Department of Mathematics} \\
{\small \it Dolivostra{\ss}e 15, 64293 Darmstadt, Germany}\\
{\small lang@mathematik.tu-darmstadt.de} \\ \\
Bernhard A. Schmitt \\
{\small \it Philipps-Universit\"at Marburg,
Department of Mathematics,}\\
{\small \it Hans-Meerwein-Stra{\ss}e 6, 35043 Marburg, Germany} \\
{\small schmitt@mathematik.uni-marburg.de}}
\maketitle

\begin{abstract}
This paper is concerned with the theory, construction and application of
variable-stepsize implicit Peer two-step methods that are super-convergent
for variable stepsizes, i.e., preserve their classical order achieved for
uniform stepsizes when applied in a gradient-based solution algorithm to
solve ODE-constrained optimal control problems in a first-discretize-then-optimize
setting. Gradients of the objective function can be computed most efficiently using 
approximate adjoint variables. High accuracy with moderate computational effort can 
be achieved through time integration methods that satisfy a sufficiently large number 
of adjoint order conditions for variable stepsizes and provide gradients with higher-order consistency.
In this paper, we enhance our previously developed variable implicit two-step Peer triplets
constructed in [J. Comput. Appl. Math. 460, 2025] to get ready for large-scale dynamical 
systems with varying time scales without losing efficiency.
A key advantage of Peer methods is their use of multiple stages with the same high stage order, 
which prevents order reduction — an issue commonly encountered in semi-discretized PDE 
problems with boundary control. Two third-order methods with four stages, good stability properties,
small error constants, and a grid adaptation by equi-distributing global errors are constructed 
and tested for a 1D boundary heat control problem and an optimal control of cytotoxic therapies 
in the treatment of prostate cancer.  
\end{abstract}


\section{Introduction}
Recently, we have developed and tested third- and fourth-order implicit Peer two-step methods \cite{LangSchmitt2022a,LangSchmitt2022b,LangSchmitt2023,LangSchmitt2024,LangSchmitt2025a}
to solve ODE-constrained optimal control problems of the form
\begin{align}
\mbox{minimize } \CC\big(y(T)\big) \label{OCprob_objfunc} &\\
\mbox{subject to } y'(t) =& \,f\big(y(t),u(t)\big),\quad
u(t)\in U_{ad},\;t\in(0,T], \label{OCprob_ODE}\\
y(0) =& \,y_0, \label{OCprob_ODEinit}
\end{align}
with the state $y(t)\in\R^m$, the control $u(t)\in\R^d$,
$f: \R^m\times\R^d\mapsto\R^m$, the objective function $C: \R^m\mapsto\R$, where the set of
admissible controls $U_{ad}\subset\R^d$ is closed and convex. The Mayer form $C(y(T))$ of the objective
function which only considers state values at $t=T$ is quite general. For example, terms given in
the Lagrange form
\begin{align}\label{objfunc-lagrange}
\CC_L(y,u) :=&\;\int_{0}^{T}l(y(t),u(t))\,dt
\end{align}
can be equivalently transferred to the Mayer form by adding a
new differential equation $y_{m+1}'(t)=l(y(t),u(t))$ and initial
values $y_{m+1}(0)=0$ to the constraints. Then \eqref{objfunc-lagrange}
simply reduces to $y_{m+1}(T)$.
\par
The design of efficient time integrators for the numerical solution of
such problems with large $m$ arising from semi-discretized time-dependent partial differential equations
is still of great interest since difficulties arise through additional adjoint order conditions, some recent literature is \cite{AlbiHertyPareschi2019,AlmuslimaniVilmart2021}.
Implicit Peer two-step methods overcome the structural disadvantages of one-step and multi-step
methods such as symplectic or generalized partitioned Runge-Kutta methods \cite{BonnansLaurentVarin2006,HertyPareschiSteffensen2013,LiuFrank2021} and backward differentiation formulas \cite{BeigelMommerWirschingBock2014}.
They avoid order reduction of one-step methods, e.g., for boundary control problems of PDEs \cite{LangSchmitt2023}, and have good stability properties.
Peer methods also allow the approximation of adjoints in a first-discretize-then-optimize (FDTO) approach with higher order,
which still seems to be an unsolved problem for multi-step methods.
FDTO is the most commonly used method and possesses the advantage of providing consistent gradients for state-of-the-art optimization algorithms. We refer the reader to the detailed discussions in our previous papers
\cite{LangSchmitt2022a,LangSchmitt2022b,LangSchmitt2023,LangSchmitt2024,LangSchmitt2025a}.
\par
So far we have considered the use of 4-stage diagonally implicit Peer two-step methods that are of order three, at least, if applied with non-uniform
time grids to a reduced first-order optimality system \cite{LangSchmitt2025a}. The transition to gradient-based solution algorithms without
analytic elimination of the control $u(t)$ requires additional order conditions for the control and certain positivity requirements for
column sums in the matrix triplet $(K_0,K,K_N)$ used in \eqref{Peerv} to define the forward Peer method \cite{LangSchmitt2024}. Starting
with two positive methods from \cite{LangSchmitt2025a}, we derive a third-order \textit{flip-symmetric pulcherrima} Peer triplet which is 
$A(61.59^o)$-stable, and a second $A(83.74^o)$-stable Peer triplet with improved stability for smoothly varying grids. Symmetry plays
an important role in the theory of geometric integration which has a strong link to the Hamiltonian structure of the first-order
optimality system discussed in Section~\ref{sec:kkt}, see \cite{HairerWannerLubich2006,SanzSerna2016}. Here, it also forms the foundation for the 
remarkably good performance of the pulcherrima triplet. Since the start and end method of the triplet are of one-step form, the requirement of high stage order prohibits the use of triangular coefficient matrices. However, for both methods  triangular iterations are provided in these steps to avoid 
additional memory requirements for large nonlinear systems.
\par
The originality of this paper lies in the development of a full-fledged, practically deployable algorithm for large-scale ODE-constrained optimal control problems. To this end, the following modifications are implemented for two academic Peer triplets from \cite{LangSchmitt2025a}: (i) improved symmetry and stability properties and (ii) memory optimization in the start and end methods. At the same time, their classical third-order accuracy is preserved even for variable step sizes.
\par
We also present a method for determining an optimized stepsize sequence \textit{a posteriori} by using error estimators to approximate global errors and implement the error equi-distribution principle introduced in \cite[Chapter 9.1.1]{AscherMattheijRussell1995}. Since Peer triplets retain their classical order even with variable time steps, these grid refinements can be applied without any loss of efficiency.
\par
The paper is organized as follows. In Section~\ref{sec:2peer}, we introduce a generalized form of a Peer two-step triplet. 
The continuous and discrete first order optimality conditions are presented in Section~\ref{sec:kkt}. The framework of gradient-based 
optimization, which relies on an iterative forward-backward marching scheme, is discussed in Section~\ref{sec:gbopt}. Several properties
of Peer triplets, including the standard and boundary methods, are described in Section~\ref{sec:propeer}. In Section~\ref{sec:2meths},
we modify two of our previous Peer two-step triplets of practical interest. Grid improvement is discussed in Section
~\ref{sec:equierr} and numerical results are given in Section~\ref{sec:num}. Conclusions are stated in Section~\ref{sec:con}.
\par
Throughout the paper, we use the following notations: $e_i$ denotes a cardinal basis vector, $\eins_k=(1,\ldots,1)\T\in\R^k$ is 
the vector of ones, and $I_k\in\R^{k\times k}$ is the identity matrix. The latter two are occasionally used without index if there is no ambiguity. Frequently used symbols are listed in Appendix~B.
\section{Two-step Peer triplets}\label{sec:2peer}
Peer two-step methods for the numerical solution of (\ref{OCprob_ODE})--(\ref{OCprob_ODEinit}), which are suitable for optimal control,
have the general redundant form
\begin{align}\label{Peerv}
 \left( A_n\otimes I_m \right) Y_n = \left( B_n\otimes I_m\right) Y_{n-1}
 + h_n \left( K_n\otimes I_m\right) F(Y_n,U_n),\ n\ge 1,
\end{align}
where $\{t_0,t_1,\ldots,t_{N+1}\}\subseteq[0,T]$ is a grid with stepsizes $h_n=t_{n+1}-t_n$, $n=0,\ldots,N$.
The stage and control solutions $Y_n=(Y_{ni})_{i=1}^s\in\R^{sm}$ and $U_n=(U_{ni})_{i=1}^s\in\R^{sd}$ are approximations
of $(y(t_n+c_ih_n))_{i=1}^s$ and $(u(t_n+c_ih_n))_{i=1}^s$, respectively, where all $s$ stages have equal accuracy.
Also, $F(Y_n,U_n)=(f(Y_{ni},U_{ni}))_{i=1}^s$ is used.
The equal stage order motivates the attribute \textit{peer} and is the key for avoiding order reduction. The nodes $c_1,\ldots,c_s$ are associated
with the interval $[0,1]$ but some may lie outside.
We note that \eqref{Peerv} is a redundant form of Peer two-step methods and may be simplified, e.g., by setting $K_n\equiv I$.
However, we will see below that the coefficient matrix $K$ leads to additional degrees of freedom when considering {\em adjoint} Peer methods.
\par
An exceptional, Runge-Kutta-like starting step is used,
\begin{align*}
\left( A_0\otimes I_m\right) Y_0=&\,a\otimes y_0 + h_0\left( K_0\otimes I_m\right) F(Y_0,U_0),
\end{align*}
with $a\in\R^s$.
With variable stepsizes, we employ as in \cite{LangSchmitt2025a} one {\em standard Peer method} $(A,K,B_n)$ for the time steps
with $1\le n< N$, where $A\in\R^{s\times s}$ is a fixed nonsingular lower triangular matrix and $K\in\R^{s\times s}$
is a fixed diagonal matrix with strictly positive entries in order to restrict the size of the stage equations to $m$ unknowns.
The matrix $B_n$ will only depend on the current stepsize ratio $\sigma_n=h_n/h_{n-1}$ in the form $B_n=B(\sigma_n)$
with a fixed function $B:\R\to\R^{s\times s}$.
In order to obtain sufficiently high orders at the boundary steps, more general matrices $A_0,A_N$ are required,
but $K_0=K_N=K$ and $B_N=B(\sigma_N)$ in the end step are also used  in this paper.
Since our focus now is more on the practical implementation of the Peer triplets, we re-designed the boundary steps
from \cite{LangSchmitt2025a}.
The increased effort for boundary steps not in diagonally-implicit form is a minor fraction of
the overall effort in the solution of the full boundary value problem only.
However, we admit that the memory requirement is multiplied there and this may lead to difficulties for large problems.
Hence we have re-designed the boundary steps and will provide fast iteration methods with triangular coefficients
$\tilde A_0,\tilde A_N$, requiring again the solution of systems with $m$ unknowns only.
This was indeed possible with positive definite diagonal $K_n\equiv K\succ0$ for all $n$.
\par
All in all, a Peer triplet now consists of a set of 5 coefficient matrices,
\[ A_0,\ A,\ A_N,\ K\succ0,\mbox{ and } B(\sigma_n),\,n=1,\ldots,N.\]
The additional matrices $\tilde A_0,\tilde A_N$ are not part of the triplet as integration method but a
supplement simplifying the solution in the boundary steps.
For a shorter notation, we will use the same symbol for a coefficient matrix like A and its Kronecker
product $A\otimes I_m$ as a mapping from the space $R^{sm}$ to itself. We also stick to $A_n$ in general
to avoid doubling of equations.

\section{First order optimality conditions}\label{sec:kkt}
The first-order optimality conditions for the optimal control problem (\ref{OCprob_objfunc})--(\ref{OCprob_ODEinit}) with
some Lagrange multiplier $p(t)\in\R^m$ read \cite{Hager2000,Troutman1996}
\begin{align}
y'(t) =& \,f\big(y(t),u(t)\big),\quad t\in(0,T],\quad y(0)=y_0, \label{KKT_state}\\
p'(t) =& \,-\nabla_y f\big(y(t),u(t)\big)\T p(t),\quad t\in[0,T),
\quad p(T)=\nabla_y {\CC}\big(y(T)\big)\T, \label{KKT_costate}\\
& \,-\nabla_u f\big(y(t),u(t)\big)\T p(t) \in N_U\big(u(t)\big),\quad t\in[0,T], \label{KKT_ctr}
\end{align}
with the normal cone mapping
\begin{align*}
N_U(u) =&\, \{ w\in\R^d: w^T(v-u)\le 0 \mbox{ for all } v\in U_{ad}\}.
\end{align*}
We recall a few facts from \cite{Hager2000}. Under appropriate regularity conditions, there exists
a local solution $(y^\star,u^\star,p^\star)$
such that the first-order optimality conditions \eqref{KKT_state}--\eqref{KKT_ctr} are necessarily satisfied.
These conditions are also sufficient, if the Hamiltonian $H(y,u,p):=p\T f(y,u)$
satisfies a coercivity assumption. Then, the control uniqueness property yields the existence
of a locally unique minimizer $u=u(\hat{y},\hat{p})$ of the Hamiltonian
over all $u\in U_{ad}$, if $(\hat{y},\hat{p})$ is sufficiently close to $(y^\star,p^\star)$.
\par
Applying a two-step Peer triplet to (\ref{OCprob_objfunc})--(\ref{OCprob_ODEinit}), we get the discrete
constraint nonlinear optimal control problem
\begin{align}
\mbox{minimize } \CC\big(y_h(T)\big) \nonumber &\\[1mm]
\mbox{subject to } A_0 Y_0=&\,a\otimes y_0+h_0KF(Y_0,U_0),\label{OCprob_peer_init}\\[1mm]
A_n Y_n=&\,B_nY_{n-1}+h_nKF(Y_n,U_n),\ n=1,\ldots,N,\label{OCprob_peer}
\end{align}
with an approximation $y_h(T)=(w\T\otimes I_m)Y_N\approx y(T)$, $w\in\R^s$, at the final time point.
Since $K=(\kappa_{ii})_{i=1}^s$ is a positive definite diagonal matrix, the first order
optimality conditions now read \cite{LangSchmitt2025a}
\begin{align}
A_0 Y_0=&\,a\otimes y_0+h_0KF(Y_0,U_0),\label{KKT_state_peer_init}\\[1mm]
A_n Y_n=&\,B_nY_{n-1}+h_nKF(Y_n,U_n),\ n=1,\ldots,N,\label{KKT_state_peer}\\[1mm]
A_N\T P_N=&\,w\otimes p_h(T)+h_NK\nabla_YF(Y_N,U_N)\T P_N,\label{KKT_adj_peer_init}\\[1mm]
A_n\T P_n=&\,B_{n+1}\T P_{n+1}+h_nK\nabla_YF(Y_n,U_n)\T P_n,\ 0\le n\le N-1,\label{KKT_adj_peer}\\[1mm]
&\,-h_n\kappa_{ii}\nabla_{u}f(Y_{ni},U_{ni})\T P_{ni}\in N_U(U_{ni}),\ 0 \le n\le N,\,i=1,\ldots,s.\label{KKT_ctr_peer}
\end{align}
Here, $p_h(T)=\nabla_y {\CC}\big(y_h(T)\big)\T$ and the Jacobians of $F$ are block diagonal matrices $\nabla_YF(Y_n,U_n)=\bldiag_i\big(\nabla_{y}f(Y_{ni},U_{ni})\big)$. Further, the adjoint solutions 
$P_n=(P_{ni})_{i=1}^s\in\R^{sm}$ are approximations of $(p(t_n+c_ih_n))_{i=1}^s$.
Since $h_n\kappa_{ii}>0$, (\ref{KKT_ctr_peer}) can be simplified to $-\nabla_{u}f(Y_{ni},U_{ni})\T P_{ni}\in N_U(U_{ni})$.
Hence, the control uniqueness property mentioned above guarantees the existence of a local minimizer $U_{ni}$ of the Hamiltonian
$H(Y_{ni},U,P_{ni})$ over all $U\in U_{ad}$, which shows the importance of having $\kappa_{ii}>0$ for all diagonal entries of $K$. Such
positivity conditions also arise in the context of classical Runge-Kutta methods or W-methods, see
\cite[Theorem~2.1]{Hager2000} and \cite[Chapter~5.2]{LangVerwer2013}.

\section{Gradient-based optimization}\label{sec:gbopt}
The cost function $\CC(y_h(T))$ can be associated with the vector of control values
\begin{align*}
U =&\,(U_{01}\T,\ldots,U_{0s}\T,U_{11}\T,\ldots,U_{Ns}\T)\T\in\R^{sd(N+1)}
\end{align*}
by interpreting $y_h(T)$ as function of $U$ and defining $\CC(U):=\CC(y_h(U))$. The whole gradient
of $\CC(U)$ can then be computed from
\begin{align*}
\nabla_{U_{ni}}C(U)=h_n\kappa_{ii}\nabla_{u}f(Y_{ni},U_{ni})\T P_{ni},
\quad 0 \le n\le N,\,i=1,\ldots,s,
\end{align*}
following, e.g., the approach from \cite{HagerRostamian1987}. Such approximate gradients
can be used to set up an iterative procedure to approximate the control vector starting
from an initial guess $U^{(0)}$. Increments $\triangle U^{(k)}$ to get updates
from
\begin{align}\label{nopt_iter}
U^{(k+1)} :=&\, U^{(k)}+\triangle U^{(k)},\quad k=0,1,\ldots
\end{align}
can be efficiently computed by gradient-based optimization algorithms such as
\textit{interior-point} \cite{ByrdGilbertNocedal2000}, \textit{sequential-quadratic-programming}
\cite{Spellucci1998}, and \textit{trust-region-reflective}
\cite{ColemanLi1994} for large-scale sparse problems with continuous objective function
and first derivatives. Various implementations are available in commercial software packages
like Matlab, Mathematica, and others.
\par
Given $U^{(k)}$, approximations for $(Y_n,P_n)$ are computable by a forward-backward marching scheme.
Since the discrete state equations (\ref{KKT_state_peer_init})--(\ref{KKT_state_peer}) do not depend on
the multipliers $P_n$, all $Y_n$ for $n=0,\ldots,N$, are computable from a simple forward calculation.
Then, using the updated values $Y_n$, one computes
$p_h(T)=\nabla_y \CC\big(y_h(T)\big)\T$ with $y_h(T)=(w\T\otimes I)Y_N$ before marching the steps 
(\ref{KKT_adj_peer_init})--(\ref{KKT_adj_peer}) backwards for $n=N,\ldots,0$, solving the discrete costate equations
for all $P_n$.
\par
Since the optimal control $u(t)$ minimizes the Hamiltonian $H(y,u,p)=p\T f(y,u)$, we may compute
an improved approximation of the control by the following minimum principle:
\begin{align}\label{ctrpp}
U^\ddagger_{ni} = \argmin_{U\in U_{ad}} H(Y_{ni},U,P_{ni}),\quad
0 \le n\le N,\ i=1,\ldots,s,
\end{align}
if $Y_{ni}$ or $P_{ni}$ are approximations of higher-order. Often, a higher convergence
order for $P_{ni}$ can be observed in practice, leading to a better approximation for
$U^\ddagger_{ni}$ compared to $U_{ni}$.

\section{Properties of Peer triplets}\label{sec:propeer}
Here, we shortly review those results from our previous paper \cite{LangSchmitt2025a} that are
necessary for the re-design of the boundary methods.
We now consider only methods with the same order $q<s$ for the original scheme and its adjoint.
\subsection{The standard method}
As for multi-step methods, order conditions are derived by Taylor expansion of the residuals if exact
solution values are used as arguments of the numerical scheme.
For variable stepsizes, these conditions depend on the stepsize ratio $\sigma_n=h_n/h_{n-1}$ in the scaling matrix $S_{n,q}=\diag_i\big(\sigma_n^{i-1}\big)\in\R^{q\times q}$.
Compact formulations of the order conditions are possible with the vector of node powers $\cc^k=(c_1^k,\ldots,c_s^k)\T$, the Vandermonde matrix $V_q=(\eins,\cc,\ldots,\cc^{q-1})\in\R^{s\times q}$, the Pascal matrix $\PP_q=\big({j-1\choose i-1}\big)\in\R^{q\times q}$, and the scaled shift matrix $\tilde E_q=\big(i\delta_{i+1,j}\big)\in\R^{q\times q}$ which commutes with $\PP_q=\exp(\tilde E_q)$.
Both the standard forward step \eqref{KKT_state_peer} and its adjoint \eqref{KKT_adj_peer} posses local order $q$ if
\begin{align}\label{OBed_vw}
 A_nV_q-KV_q\tilde E_q=&B_nV_q\PP_q^{-1}S_{n,q}^{-1}\equiv B(1)V_q\PP_q^{-1},\\\label{OBed_ad}
 A_n\T V_q+KV_q\tilde E_q=&B_{n+1}\T V_qS_{n+1,q}\PP_q\equiv B(1)\T V_q\PP_q.
\end{align}
Since the matrices on the left-hand sides of both equations are independent of $\sigma_n,\sigma_{n+1}$ by design, the same holds for the other terms which is highlighted by the right-most expressions.
Combining both conditions with varying parameter $\sigma$ leads to a very restricted form of the following matrix playing a central role \cite[Lemma~3.1]{LangSchmitt2025a},
\begin{align}\label{MatQ}
 \QQ_{q,q}(\sigma)=e_1e_1\T\in\R^{q\times q},
 \quad\mbox{where } \QQ_{q,r}(\sigma):=V_q\T B(\sigma)V_r\PP_r^{-1},\ q,r\le s.
\end{align}
This property is also the reason why the uniform local order should be smaller than $s$, since for $q=r=s$ no freedom is left for fulfilling $\sigma$-dependent conditions.
As a consequence it turned out that the following congruent $s\times s$ - matrices
\begin{align}\label{kongruent}
 \hat A_n=V_s\T A_nV_s,\quad \hat K:=V_s\T KV_s,\quad \hat B(\sigma):=V_s\T B(\sigma) V_s
\end{align}
also possess very restricted forms with few free parameters.
Since we are mainly interested here in four-stage methods, we depict this form for the standard method with order $q=3=s-1$:
\begin{align}\label{AdBd}
 \hat A=\begin{pmatrix}
  1&1&1&1\\
  0&\frac12&\frac23&\hat a_{24}\\[1mm]
  0&\frac13&\frac24&\hat a_{34}\\[1mm]
  \hat a_{41}&\hat a_{42}&\hat a_{43}&\hat a_{44}
 \end{pmatrix},\quad
  \hat B(\sigma)=\begin{pmatrix}
  1&1&1&1\\
  0&0&0&\hat b_{24}(\sigma)\\[1mm]
  0&0&0&\hat b_{34}(\sigma)\\[1mm]
  \hat a_{41}&\hat b_{42}(\sigma)&\hat b_{43}(\sigma)&\hat b_{44}(\sigma)
 \end{pmatrix},
\end{align}
and $\hat K=hankel(1,\frac12,\frac13,\frac14,\hat k_5,\hat k_6,\hat k_7)$ is a Hankel matrix coinciding
with the Hilbert matrix in its first 4 antidiagonals.
The first row and column of $\hat B(\sigma)$ coincide with that of $\hat A$ due to the
conditions \eqref{OBed_vw}, \eqref{OBed_ad} of lowest order, and with the exception of $\hat b_{44}(\sigma)$
its other elements are fixed by the order $q=3$.
Since both methods in this paper will also satisfy \eqref{LSRKv} below, we have already fixed $\hat a_{14}=1$.
A drawback of the sparing parametrization through \eqref{kongruent}, \eqref{AdBd} is that the triangular form
of $A$ itself is lost and has to be enforced by an appropriate set of equations.
As a consequence the only free parameters of the standard Peer method will be $\hat a_{41}$ and $\hat b_{44}(\sigma)$
in \eqref{AdBd} and the nodes $c_i$, of course.
\par
A critical requirement for the standard method is zero stability which may be very demanding for variable stepsizes.
For the trivial ODE $y'=0$, the stability matrix of the forward step \eqref{KKT_state_peer}, $Y_n=A^{-1}B_nY_{n-1}$, of
the standard method is $\bar B_n:=A^{-1}B_n$ and for the adjoint step \eqref{KKT_adj_peer} it is 
$\tilde B_{n+1}\T:=A\mT B_{n+1}\T=A\mT\bar B_{n+1}\T A\T$.
Zero stability requires that arbitrarily long products, e.g., $\bar B_n\cdots \bar B_{k+1},\,n>k,$ are uniformly bounded.
Since $\hat A^{-1}\hat B(\sigma)=V_s^{-1}\bar B(\sigma) V_s$ is similar to the stability matrix $\bar B(\sigma)$,
based on the representation \eqref{AdBd} in \cite{LangSchmitt2025a}, the construction of a $\sigma$-independent
weight matrix $W\in\R^{s\times s}$ was incorporated into the design of all standard methods, which block-diagonalizes all matrices simultaneously,
\begin{align}\label{decoup}
 W^{-1}\bar B(\sigma)W=\begin{pmatrix}
 1&0\\0&B_{se}(\sigma)
 \end{pmatrix},\mbox{ such that }\|B_{se}(\sigma)\|_\infty\le\tilde\gamma<1
\end{align}
for $\sigma$ near one.
In particular, \eqref{decoup} states that 1 is an eigenvalue of $\bar B(\sigma)$.
This follows from the order conditions \eqref{OBed_vw}, \eqref{OBed_ad}, from which the exact form of the
corresponding eigenvectors can be deduced as
\begin{align}\label{eigenv}
 \bar B(\sigma)\eins=A^{-1}B(\sigma)\eins=\eins,\quad
 (\eins\T A) \bar B(\sigma)=(\eins\T A).
\end{align}
The weighted norm corresponding to \eqref{decoup} satisfies
\begin{align}\label{nullstab}
  \ltnorm \bar B(\sigma)\rtnorm:=\|W^{-1}\bar B(\sigma)W\|_\infty=1,\ \sigma\in[\underline{\sigma},\bar\sigma],\ 0<\underline\sigma<1<\bar\sigma,
\end{align}
where a large interval $[\underline{\sigma},\bar\sigma]$ was one of the design objectives.
Based on this property, convergence of the global error with order $q-1=2$  was proven rigorously in
\cite[Theorem~6.1]{LangSchmitt2025a}.
\par
In order the regain global order $q=3$, an additional super-convergence property is employed, which prohibits the propagation of the leading term of the local error.
Denoting evaluations of the exact solution by boldface as in ${\bf y}_n=\big(y^\star(t_{ni})\big)_{i=1}^n$, ${\bf y}'_n=\big((y^\star)'(t_{ni})\big)_{i=1}^n$, etc., the  local errors of the standard method are defined by
\begin{align}\label{locerry}
 \tau_n^Y=&\;{\bf y}_n-A_n^{-1}\big(B_n{\bf y}_{n-1}+h_nK{\bf y}_n'\big),\; n=1,\ldots,N,\\\label{locerrp}
 \tau_n^P:=&\;{\bf p}_n-A_n\mT\big(B_{n+1}\T{\bf p}_{n+1}-h_nK{\bf p}_n'\big),\; n=0,\ldots,N-1,
\end{align}
and the order conditions \eqref{OBed_vw}, \eqref{OBed_ad} lead to the following expressions of these errors \cite{LangSchmitt2025a}:
\begin{align}\label{tauyn}
 \tau_n^Y=&\;h_n^q\beta_{q,n}(\sigma_n)(y^\star)^{(q)}(t_n)+O\left( h_n^{q+1}\|(y^\star)^{(q+1)}\|_{[n]}\right),\\\label{betav}
  &\beta_{q,n}(\sigma):=\frac1{q!}A_n^{-1}\big(A_n\cc^q-B_n(\cc-\eins)^q\sigma^{-q}-qK\cc^{q-1}),\ n\ge 1,\\\label{taupn}
  \tau_n^P=&\;h_n^q\beta_{q,n}^\dagger(\sigma_{n+1})(p^\star)^{(q)}(t_n)+O\left( h_n^{q+1}\|(p^\star)^{ (q+1)}\|_{[n]}\right),\\\label{betaa}
  &\beta_{q,n}^\dagger(\sigma)=\frac1{q!}A_n\mT(A_n\T\cc^q-B_{n+1}\T(\eins+\sigma\cc)^q+qK\T\cc^{q-1}),\,n<N.
\end{align}
As the theoretical basis of stepsize control, we emphasize here the localized dependence on derivatives and write the interval $[0,T]=\bigcup_{n=0}^N\pi_n$ as the union of the subintervals $\pi_n:=[t_n,t_{n+1}]$.
To take into account the two-step form of the Peer methods, for $\phi\in C[0,T]$, we use the definition
\begin{align}\label{nachbarn}
 \|\phi\|_{[n]}:=\max\{\|\phi(t)\|:\,t\in\pi_{n-1}\cup\pi_n\cup\pi_{n+1}\},\ 0\le n\le N,
\end{align}
where $\pi_{-1}=\pi_{N+1}=\emptyset$.
Through the iterated recursion \eqref{KKT_state_peer}, the local error $\tau_k^Y$ is multiplied by long products $\bar B_n\cdots\bar B_{k+1}\tau_k^Y$ which by \eqref{decoup}, \eqref{eigenv} converge as
\[ \bar B_n\cdots\bar B_{k+1}=\eins\eins\T A+O(\tilde\gamma^{n-k}),\ n-k\to\infty,\]
see \eqref{decoup}, with an analogous version for the adjoint matrices $\tilde B_{n+1}\T$.
Hence, there is no accumulation of leading error terms of the standard method under the super-convergence conditions
\begin{align}\label{supkonv}
\eins\T A\beta_q(\sigma)=0,\quad \eins\T A\T \beta_q^\dagger(\sigma)=0,
\end{align}
and the global order $q=s-1=3$ is regained, see \cite[Lemma~6.1]{LangSchmitt2025a}.
By purpose, we did not specify here the range for the stepsize ratio $\sigma$.
In fact, if \eqref{supkonv} is satisfied for arbitrary $\sigma\in\R$, the standard Peer method (as in \texttt{AP4o33vgi} below)
will converge for \textit{general grids} satisfying only the restriction in \eqref{nullstab}.
On the other hand, if \eqref{supkonv} is satisfied for $\sigma=1$ only, order $q=s-1=3$ still holds for \textit{smooth grids} satisfying
\begin{align}\label{glatt}
 |\sigma_n-1|\le\eta{h_n},\ n=0,\ldots,N,
\end{align}
with some moderately large constant $\eta$ (e.g. $\eta=15$). The other triplet \texttt{AP4o33vsi} below belongs to this class.
As a trade-off for its inferior convergence properties, it possesses a much larger angle of $A(\alpha)$-stability than \texttt{AP4o33vgi}.
\par
In the design process in \cite{LangSchmitt2025a} it turned out that a special feature of the standard method leads to improved properties concerning accuracy and stability. In addition to the order conditions \eqref{OBed_vw}, \eqref{OBed_ad} with $q=s-1$, these methods satisfy
\begin{align}\label{LSRKv}
 e_1\T\hat A=\eins\T,\quad c_s=1.
\end{align}
An essential consequence is that $\eins\T A=e_s\T$ holds, leading to the simple left eigenvector $e_s\T\bar B(\sigma)=e_s\T$ of the stability matrix.
Hence, the last stage of the Peer method may also be written as
\begin{align}\label{LSRK}
Y_{ns}=Y_{n-1,s}+h_n\sum_{j=1}^s\kappa_{jj}f(Y_{nj},U_{nj}),
\end{align}
which, together with $c_s=1$, resembles a final stage of a Runge-Kutta method since $Y_{n-1,s}\cong y(t_n)$.
We say that such methods have LSRK form ({\em Last Stage is Runge-Kutta}).
However, we stress the fact that this last Runge-Kutta step is not prone to order reduction since the stage
increments $f(Y_{nj},U_{nj})$ have full stage order $q$ due to the two-step structure of the Peer method.
\subsection{Boundary methods}
Exceptional coefficients $A_0,A_N$ in general form are required for the boundary steps since with triangular matrices the first starting stage in \eqref{KKT_state_peer_init} and the first backward stage (with index $s$) of the adjoint end condition \eqref{KKT_adj_peer_init} would be implicit Euler steps with insufficient local order two.
Since the boundary steps are applied only once each, their required local order is identical with the global order $q$.
For the starting step, the forward order condition reads \cite{LangSchmitt2025a}
\begin{align}\label{OBedstrt}
 A_0V_q=ae_1\T+KV_q\tilde E_q.
\end{align}
Since $\tilde E_qe_1=0$, we deduce that $a=A_0\eins$.
The adjoint order condition \eqref{OBed_ad} for the starting step may be compared to the one of the standard method, yielding
\begin{align*}
 A_0\T V_q=B(\sigma_1)\T V_qS_q(\sigma_1)\PP_q-KV_q\tilde E_q=A\T V_q.
\end{align*}
According to \eqref{OBed_ad}, here all 4 product terms are $\sigma$-independent, showing that $V_q\T(A_0-A)=0$.
Hence, for $q=s-1$, the matrix $A_0$ is a rank-1-modification of $A$ itself,
\begin{align*}
 A_0=A+ V_s\mT e_s\phi_0\T,\quad \phi_0\in\R^s,
\end{align*}
since the last row of $V_s^{-1}$ is orthogonal to all columns of $V_q,\,q<s$.
The coefficient $a$ now is $a=A\eins+(e_1\T\phi_0)V_s\mT e_s$ and with respect to \eqref{OBed_vw},
condition \eqref{OBedstrt} may be rewritten as
\begin{align*}
(A\eins+(e_1\T\phi_0)V_s\mT e_s)e_1\T-V_s\mT e_s\phi_0\T=AV_q-KV\tilde E_q=B(\sigma_1)V_q\PP_q^{-1}S_q(\sigma_1)
\equiv B(1)V_q\PP_q^{-1}.
\end{align*}
Multiplying by $V_s\T$ from the left and recalling that $V_s\T A\eins=\hat Ae_1$ gives
\begin{align}\label{StrtNB}
 \hat A e_1e_1\T+e_s\big((e_1\T\phi_0)e_1-\phi_0\big)\T =V_s\T B(1)V_q\PP_q^{-1}=\QQ_{s,q}(1),
\end{align}
with the matrix $\QQ$ from \eqref{MatQ}.
Since $\QQ_{s,q}e_1=\hat Be_1=\hat Ae_1$, the first column in \eqref{StrtNB} is satisfied for arbitrary $\phi_0$ and the same holds for the first $q$ rows, since $\QQ_{q,q}=e_1e_1\T$ by \eqref{MatQ}.
The remaining entries yield $e_j\T\phi_0=-e_s\T\QQ_{s,q}e_j$, $j=2,\ldots,q$.
\par
In an analogous way, the forward condition for the end step \eqref{KKT_state_peer} means that
\begin{align*}
 A_NV_q=B_NV_q\PP_q^{-1}S_{N,q}^{-1}-KV_q\tilde E_q=AV_q,
\end{align*}
leading to $(A_N-A)V_q=0$, and for $q=s-1$ to
\begin{align*}
 A_N=A+\phi_Ne_s\T V_s^{-1}
\end{align*}
with parameter $\phi_N\in\R^s$.
The order condition for the adjoint starting step \eqref{KKT_adj_peer_init} is \cite{LangSchmitt2025a}
\begin{align}\label{OBedend}
 A_N\T V_q+KV_q\tilde E_q=w\eins_q\T,
\end{align}
revealing that $w=A_N\T V_se_1=A_N\T\eins=A\T\eins+V_s\mT e_s\phi_N\T\eins$.
Multiplying \eqref{OBedend} by $V_s\T$ from the left and reminding \eqref{OBed_ad}, we obtain
\begin{align}\notag
  V_s\T w\eins_q\T=&V_s\T(A\T\eins+V_s\mT e_s\phi_N\T\eins)\eins_q\T=
 \\[1mm]\label{EndNB}
 =\hat A\T e_1\eins_q\T+(\phi_N\T\eins)e_s\eins_q\T
 \stackrel{\eqref{OBedend}}=&V_s\T(A\T V_q+KV_q\tilde E_q)+e_s\phi_N\T V_q
 =V_s\T B(1)\T V_q\PP_q+e_s\phi_N\T V_q.
\end{align}
Multiplying by $\PP_q^{-1}$ from the right and reminding $\eins\T\PP_q^{-1}=e_1\T$, the first $q=s-1$ rows of the lower equation have the form $\eins e_1\T\in\R^{q\times q}$ on both sides by \eqref{AdBd}.
And the first entry of row number $s$ is satisfied for arbitrary $\phi_N$, but the remaining entries lead to $q-1$ restrictions.
We collect these results now.
\begin{lemma}\label{LRng1}
For $q=s-1$ let $A_0$ satisfy the order conditions \eqref{OBedstrt}, \eqref{OBed_ad} for $n=0$ and $A_N$ the conditions \eqref{OBedend}, \eqref{OBed_vw} for $n=N$ and let \eqref{MatQ} hold.
Then, with $\phi_0,\phi_N\in\R^s$ it holds
\begin{align}\label{Rdrng}
 A_0=A+V_s\mT e_s\phi_0\T,\quad A_N=A+\phi_Ne_s\T V_s^{-1},
\end{align}
where
\begin{align}\label{NBphi}
 e_j\T \phi_0=-e_s\T\QQ_{s,q}(1)e_j,\quad
 \phi_N\T V_q \PP_q^{-1}e_j =-\hat b_{js}(1),\quad j=2,\ldots,q.
\end{align}
\end{lemma}
The Lemma shows that in the boundary methods there are only 2 free parameters in the rank-one representation \eqref{Rdrng} of $A_0,A_N$, each.
\begin{remark}
The fact that the conditions for $\phi_0,\phi_N$ can be fulfilled although the left hand sides of \eqref{StrtNB} and \eqref{EndNB} have rank 2 only while the matrix $B(1)$ on the right hand side may be a full matrix, in principle, may seem like a lucky coincidence.
However, the analysis in \cite{LangSchmitt2024} revealed the strong relation between the coefficients of the standard and the boundary methods forced by the many order conditions.
The design of the standard methods has taken account for these side conditions.
\end{remark}
Finally, we have to complement the definition of the local errors \eqref{tauyn}, \eqref{taupn} for the exceptional order conditions \eqref{OBedstrt}, \eqref{OBedend} at the boundaries through
\begin{align}\label{betarand}
 \beta_{q,0}=\frac1{q!}\big(\cc^q-qA_0^{-1}K\cc^{q-1}\big),\quad
 \beta_{q,N}^\dagger=\frac1{q!}\big(\cc^q+qA_N\mT K\cc^{q-1}-\eins\big).
\end{align}
These coefficients do not depend on $\sigma$.
\subsection{Triangular iterations at the boundaries}
Although the two boundary steps require a fraction of the overall computational effort only, with full matrices $A_0,A_N$ the memory requirement for nonlinear systems of size $(sm)\times(sm)$ may be prohibitive for large-scale problems.
This difficulty may be eliminated by reviving an idea from our paper \cite{LangSchmitt2022a} employing iterations with triangular approximations $\tilde A_0,\tilde A_N$ of the matrices $A_0,A_N$.
Implementation of this iteration requires only minor changes to the standard step if the subdiagonals of $\tilde A_0$ and $A_0$ as well as of $\tilde A_N$ and $A_N$ coincide, which means that
\begin{align}\label{Randapp}
 R_0:=A_0-\tilde A_0,\quad R_N:=A_N-\tilde A_N
\end{align}
are upper triangular.
In this case, the boundary steps can be solved stage by stage and approxima\-tions for previous stages can be overwritten immediately, avoiding additional memory requirements.
We shortly describe this property for an approximate Newton step of the starting method \eqref{KKT_state_peer_init}. Here,
\begin{align}\label{Nwtit}
 \big(\tilde A_0-h_0 K\textbf{J}_0\big)(Y_0^{(k+1)}-Y_0^{(k)})=a\otimes y_0+h_0KF(Y_0^{(k)},U_0)-A_0Y_0^{(k)},\ k\ge 0,
\end{align}
with the Jacobian $\textbf{J}_0=\bldiag_i(J_{0i})$, $J_{0i}=\nabla_y f(y_0,e_1\T V_4^{-1}U_0),\ i=1,\ldots,s$.
Since $\tilde A_0=\big(\tilde a_{ij}^{(0)}\big)$ is lower triangular, the equation for $Y_{0i}^{(k+1)},\,i=1,\ldots,s,$ may be solved successively as
\begin{align}\notag
 &\big(\tilde a_{ii}^{(0)}I-h_0\kappa_{ii}J_{0i}\big)(Y_{0i}^{(k+1)}-Y_{0i}^{(k)})\\\notag
 &=-\sum_{j=1}^{i-1}\tilde a_{ij}^{(0)}(Y_{0j}^{(k+1)}-Y_{0j}^{(k)})
 +a_iy_0+h_0\kappa_{ii}f(Y_{0i},U_{0i})-\sum_{j=1}^sa_{ij}^{(0)}Y_{0j}^{(k)}\\\notag
 &=a_iy_0+h_0\kappa_{ii}f(Y_{0i},U_{0i})-\sum_{j=1}^{i-1}\tilde a_{ij}^{(0)}Y_{0j}^{(k+1)}-\sum_{j=i}^{s} 
 \left( r_{ij}^{(0)} +\delta_{ij}\tilde a_{ii}^{(0)} \right) Y_{0j}^{(k)}.
\end{align}
Obviously, the old value $Y_{0i}^{(k)}$ may be overwritten by $Y_{0i}^{(k+1)}$ after solving stage equation number $i$.
\par
For linear autonomous problems, the error of the iteration for the end equations obeys the recursion
$\check Y_n^{(k+1)}=(\tilde A_n-h_nK\otimes J_n)^{-1}(\tilde A_n-A_n)\check Y_n^{(k)}$, $n\in\{0,N\}$.
Considering an eigenvalue $\lambda$ of $J_n$ and $z=h_n\lambda$, the linear convergence of the iteration is described
by the matrix
\begin{align}\label{ItMat}
 \Sc_n(z):=(\tilde A_n-zK)^{-1}(\tilde A_n-A_n),\ z\in\C.
\end{align}
For the two Peer triplets discussed below, new boundary methods have been constructed by minimizing the spectral radius $\rho_\R$ of the iteration matrix \eqref{ItMat} along the negative real axis, resulting in small contraction factors below $0.1$.
Also the spectral radius $\rho_\alpha$ in the sector of $A(\alpha)$-stability will be reported, since this sector is the largest one where an application of the triplet may be advised.
Since these contractions factors are quite small, they might scarcely slow down the convergence of the Newton iteration.
\par
These contraction factors are also valid for the adjoint iterations.
The iteration matrix for \eqref{KKT_adj_peer}, for instance, is
\[ \Sc^\dagger_0(z)=(\tilde A_0\T-zK)^{-1}(\tilde A_0-A)\T\]
and possesses the same eigenvalues as $\Sc_0(z)$ since eigenvalues are invariant if two matrix factors are commuted.
\begin{remark}
Replacing the matrix $A_0-h_0 K\textbf{J}_0$ in the simplified Newton method by the lower block-triangular
matrix $\tilde A_0-h_0 K\textbf{J}_0$ with a small contraction factor for linear autonomous problems still
delivers a fast converging sequence of approximate solutions. Thus, with appropriate stopping criteria, the 
overall third order of the Peer triplets
is not affected by the triangular approximations at the boundaries.
\end{remark}

\section{Two modified Peer triplets}\label{sec:2meths}
Here we present modifications of two Peer triplets from \cite{LangSchmitt2025a} with definite matrices $K\succ0$.
In order to avoid ambiguities, the letter \textbf{i} (for iteration in the boundary methods) is added to their old names.
Beyond the order of convergence, important properties of the triplets are $A(\alpha)$-stability and zero-stability \eqref{nullstab},
the leading error constants
\[ err_{q,n}=\|\beta_{q,n}(1)\|_\infty,\quad
  err_{q,n}^\dagger=\|\beta_{q,n}^\dagger(1)\|_\infty,
\]
and the precise form of the super-convergence condition \eqref{supkvvg}.
Also, the Jacobians $A_n-h_nK{\bf J}_n=K(K^{-1}A_n-h_n{\bf J}_n)$ of the original boundary steps should be well-conditioned for stiff problems with eigenvalues of ${\bf J}_n$ in the left complex half-plane, and we require the eigenvalue conditions
\begin{align}
 \mu_n:=\min_j\text{Re}\,\lambda_j(K^{-1}A_n)>0,\ n\in\{0,N\}.
\end{align}

\subsection{The 'pulcherrima' triplet \texttt{AP4o33vgi}}
For the triplet \texttt{AP4o33vgi} ({\textbf A}djoint {\bf P}eer method with \textbf{4} stages, \textbf{o}rders $q=\textbf{3}$ for state and adjoint, \textbf{v}ariable stepsizes, \textbf{g}eneral grids, \textbf{i}teration), the standard method is given with the nodes and coefficients
\begin{align}\label{pulchcka}
 \cc\T=\left(0,\frac13,\frac23,1\right),\quad
 K=\diag\left(\frac18,\frac38,\frac38,\frac18\right)\succ0,\quad
 A=\begin{pmatrix}
  1&\cdot&\cdot&\cdot\\[1mm]
  -\frac94&\frac94&\cdot&\cdot\\[1mm]
  \frac94&-\frac92&\frac94&\cdot\\[1mm]
  -1&\frac94&-\frac94&1
 \end{pmatrix}.
\end{align}
Vanishing entries are replaced by dots in order to highlight the structure.
The nodes and the diagonals of $K$ belong to the pulcherrima quadrature rule having positive weights.
Since $c_4=1$ and $\eins\T A=e_s\T$, it has LSRK form.
A further remarkable property is that the standard method is its own adjoint, which means that if the grid on the interval
$[0,T]$ is flipped at its midpoint, the original scheme is reproduced.
Since this property is of interest, we require it for the whole triplet.
Algebraically, it may be described with the flip permutation $\Pi:=\big(\delta_{i,s+1-i}\big)\in\R^{s\times s}$.
The nodes are symmetric, $\Pi\cc=\eins-\cc$ and the coefficients satisfy
\begin{align}\label{selfadj}
 \Pi A\Pi=A\T,\; \Pi K\Pi=K,\; \Pi B(\sigma)\Pi=B(\sigma^{-1})\T, \;
 w=\Pi a,\; \Pi A_N\Pi=A_0\T.
\end{align}
The condition for $B(\sigma)$ is more general than in \cite{LangSchmitt2025a}, where it was used for $\sigma=1$ only and is based on the fact that by a flip of the grid stepsize ratios are inverted.
Hence, the triplet \texttt{AP4o33vgi} also has a slight modification in the last element of the third coefficient of the standard method
\begin{align}\label{b44sym}
\hat b_{44}(\sigma)=&\frac1{804}\left(132\sigma+\frac{65}{\sigma}-149\right),\\[2mm]
 \hat B(\sigma)=&\begin{pmatrix}
  1&1&1&1\\[1mm]\label{Bsym}
  \cdot&\cdot&\cdot&\frac1{36\sigma}\\[1mm]
  \cdot&\cdot&\cdot&\cdot&\\[1mm]
  0&\frac{\sigma}{36}&\frac{\sigma}{18}&\hat b_{44}(\sigma)
 \end{pmatrix},
 \quad
 W=\begin{pmatrix}
 1& -2& \frac{24}5& -\frac92\\[1mm]
 1& -\frac43& 0& \frac32\\[1mm]
 1& -\frac23& -\frac85 & \frac32\\[1mm]
 1& 0& 0& 0
 \end{pmatrix}.
\end{align}
The new coefficient also leads to a slightly larger interval $[\underline{\sigma},\bar\sigma]=[0.57,2.1]$ of uniform zero-stability \eqref{nullstab} and equation \eqref{Bsym} contains a better conditioned version of the corresponding weight matrix $W$.
The stability region $\{z\in\C:\,\varrho\big(M(z)\big)\le1\},\;M(z):=(A-zK)^{-1}B,$ has a simple form.
Its complement is kidney-shaped, contained in $[-4,31]\times[-22,22]$, and touches the imaginary axis at $z=0$ due to consistency.
The method is $A(\alpha)$-stable with $\alpha=61.59^o$.
A remarkable property is that it even satisfies $\varrho(M(i\xi))\le1$ for $\xi\in[-1,1]$.
\par
Since $c_1=0$ and $\hat Ae_1=e_1$, this triplet possesses LSRK form also for the adjoint equation, but this is also clear from the flip condition \eqref{selfadj}.
Important consequences of the LSRK property are that super-convergence occurs for arbitrary stepsize ratios $\sigma$:
\begin{align}\label{supkvvg}
 \eins\T A\beta_3(\sigma)=e_4\T\beta_3(\sigma)\equiv 0,\quad
 \eins\T A\T \beta_3^\dagger(\sigma)=e_1\T\beta_3^\dagger(\sigma)\equiv0,\quad \sigma\in\R,
\end{align}
for the leading error coefficients \eqref{betav}, \eqref{betaa}.
Essential properties of the standard method, the stage number and order, range of nodes, stability angle $\alpha$, interval of uniform zero stability and error constants, are collected in Table~\ref{TPT}.
\begin{table}
\centerline{\begin{tabular}{|l|c|c|c|c|c|c|c|c|}\hline
  triplet &$s,\,q$ & nodes & $\alpha$ & $\sigma\in$&$err_{3}$&$err_3^\dagger$\\\hline
  AP4o33vgi&$4,\,3$&$[0,1]$&$61.59^o$&$[0.57,2.10]$&9.8e-3&9.8e-3\\
  AP4o33vsi&$4,\,3$&$(0,1]$&$83.74^o$&$[0.65,1.80]$&5.1e-2&3.2e-2\\
\hline
\end{tabular}}
\caption{Properties of the standard methods.
}\label{TPT}
\end{table}
\begin{table}
\centerline{\begin{tabular}{|l|c|c|c|c|c|c|c|c|c|c|}\hline
 &\multicolumn{5}{c|}{Starting method}&\multicolumn{5}{c|}{End method}\\\cline{2-11}
  triplet &$\rho_{\R,0}$&$\rho_{\alpha,0}$& $\mu_0$&$err_{3,0}$& $err_{3,0}^\dagger$& $\rho_{\R,N}$ & $\rho_{\alpha,N}$ & $\mu_N$ & $err_{3,N}$&$err_{3,N}^\dagger$\\\hline
  AP4o33vgi& 6.4e-2 & 0.155& 4.31 &5.2e-3 &9.5e-3
  &6.4e-2 & 0.155  & 4.31&9.5e-3 &5.2e-3\\
  AP4o33vsi&3.4e-2 &0.126 & 5.65 &5.2e-3 & 2.1e-2
   & 6.6e-2& 0.217&2.55& 6.7e-2&4.1e-2\\\hline
\end{tabular}}
\caption{Properties of the boundary methods.
}\label{TPB}
\end{table}
\par
The boundary methods are re-designed now in order to obtain fast iterations \eqref{Nwtit} with triangular approximations $\tilde A_0,\tilde A_N$.
Not only the free diagonal parameters of these matrices will be used, but also the two free parameters of $A_0$ resp. $A_N$ itself.
For the rank-one-representations \eqref{Rdrng}, the flip condition \eqref{selfadj} means that
\begin{align*}
 \Pi A_N\Pi =\Pi A\Pi +\Pi\phi_N e_s\T V_s^{-1} \Pi\stackrel!=A_0\T=A\T+\phi_0e_s\T V_s^{-1} \Pi\folgt \phi_0=\Pi\phi_N,
\end{align*}
and for the starting method, we simply use $\tilde A_0=\Pi\tilde A_N\Pi$.
As mentioned before, for the triangular iteration, we consider the contraction factors
\begin{align}\label{kntrkt}
 \rho_{\R,n}:=\max\{\varrho\big(S_n(z)\big):\ z\in(-\infty,0)\},\quad
 \rho_{\alpha,n}:=\max\{\varrho\big(S_n(z)\big):\ |\arg(z)-\pi|\le\alpha\},
\end{align}
where $\varrho(\cdot)$ means the spectral radius and $\alpha$ is the angle of $A(\alpha)$-stability of the method listed in Table~\ref{TPT}.
In fact, only the contraction factor $\rho_{\R,N}$ on the negative real axis was minimized but the factor $\rho_{\alpha,N}$ still is acceptable.
The matrices for \texttt{AP4o33vgi} are printed in Appendix~A1 and the contraction factors are
 $\rho_{\R,0}=\rho_{\R,N}=0.0637,\quad\rho_{\alpha,0}=\rho_{\R,N}=0.155$, see also Table~\ref{TPB}.
\par
Since $c_s=1$, the best approximation of order $s$ at the end point would be obtained with $w\T=e_s\T=\eins\T V_s^{-1}$ by $y_h(T)=e_s\T Y_N=Y_{Ns}$.
However, \eqref{OBedend} requires that $w\T=\eins\T A_N$ or
$$w\T V_s=\eins\T\hat A_N V_s\stackrel{\eqref{Rdrng}}=e_1\hat A+\eins\T\phi_N e_s\T=(1,\ldots,1,1+\eins\T\phi_N).
$$
This means that $w\T Y_N=e_s\T Y_{N}+(\eins\T\phi_N)e_s\T V_s^{-1} Y_N=Y_{Ns}+O(h_N^{s-1})$ is an $O(h_N^q)$ perturbation of the 'best' approximation with the order $q$ of the triplet and the choice $w=e_s$ would give away one useful degree of freedom in the end method.
\subsection{The triplet \texttt{AP4o33vsi} for smooth grids}
Although the pulcherrima triplet has many favourable properties, its stability angle below $62^o$ may be too small for some applications.
By sacrificing some of these properties, a Peer triplet with much larger angle could be designed in \cite{LangSchmitt2025a}.
Being designed with algebraic manipulation software, all coefficients of this triplet are known with exact rational numbers.
However, we mostly present double precision data in order to save space.
The triplet \texttt{AP4o33vsi} ({\textbf A}djoint {\bf P}eer method with \textbf{4} stages, \textbf{o}rders $q=\textbf{3}$ for state and adjoint, \textbf{v}ariable stepsizes, \textbf{s}mooth grids, \textbf{i}teration) is based on the nodes
\begin{align}\label{Knotn3vs}
 \cc\T=\left(\frac{144997}{389708},\frac{73}{748},\frac{77297572}{117896267},1\right)
 \doteq(0.37207,0.09759,0.65564,1).
\end{align}
The diagonal matrix $K$ is also positive definite and, hence, preserves the sign in the cone condition \eqref{KKT_ctr_peer}.
The triplet has a substantially larger
angle $\alpha>83^o$ of $A(\alpha)$-stability than \texttt{AP4o33vgi} and the complement of its stability region is contained in $[-0.6,20]\times[-15,15]$.
One drawback is larger error constants, see Table~\ref{TPT}.
The standard method is uniformly zero stable in $[\underline\sigma,\bar\sigma]=[0.65,1.80]$ 
with the modified weight matrix
\[ W=\begin{pmatrix}
  1& -\frac7{18}& -\frac2{21}& \frac9{14}\\[1mm]
  1& -\frac{47}{84}& \frac{11}{14}& -\frac{95}{126}\\[1mm]
  1& -\frac3{14}& -\frac{26}{63}& \frac{17}{30}\\[1mm]
  1& 0& 0& 0
\end{pmatrix}.\]
The method has LSRK form again, but not its adjoint.
For the adjoint method, we now have in condition \eqref{supkonv} for super-convergence \cite[Lemma~5.1]{LangSchmitt2025a},
\[ \eins\T A\T\beta_q^\dagger(\sigma)=\hat a_{s1}(1-\sigma^q),\]
where $\hat a_{41}\le0.1011$ here.
Hence, for \texttt{AP4o33vsi} it only holds that
\begin{align}\label{supkvvs}
 \eins\T A\beta_3(\sigma)=e_4\T\beta_3(\sigma)\equiv 0,\quad
 |\eins\T A\T\beta_3^\dagger(\sigma_n)|=|\hat a_{41}(1-\sigma^3)|
 \le\frac23|1-\sigma|
 \ \le\frac23\eta h_n,
\end{align}
if $\sigma_n\in[\underline\sigma,\bar\sigma]$ and \eqref{glatt} are obeyed in the grid construction.
The constant in \eqref{supkvvs} is computetd from $\hat a_{41}(1+\bar\sigma+\bar\sigma^2)\doteq0.611\le2/3$.
This means that the adjoint method possesses global order $q=3$ for smooth grids only.
\par
Since \texttt{AP4o33vsi} is not flip-symmetric, the boundary methods are independent.
For both methods, optimization of the contraction factors was performed for each with the 4 free diagonals of $\tilde A_0$ resp. $\tilde A_N$ and the 2 additional parameters from the boundary methods itself, see \eqref{Rdrng}.
The coefficients of these methods are given in Appendix~A2.
The contraction factors $\rho_{\R}$ on the real axis are again below $1/15$ and despite the much larger value of the stability angle $\alpha>83^o$ also $\rho_{\alpha,N}\doteq 0.22$ may be acceptable while $\rho_{\alpha,0}\cong1/8$ is even smaller than for the pulcherrima triplet.
Also the error constants of the boundary methods have a size similar to the standard method, see Table~\ref{TPB}.
\subsection{Global error estimates}\label{ssec:globerr}
The global error in the unconstrained case $N_U=\{0\}$ may be derived by a slight generalization of the proof presented in \cite{LangSchmitt2024} for uniform grids to the case of general grids obeying certain restrictions.
The key point is that the design of the Peer methods includes the construction of fixed weight matrices $W,W^\dagger=(AW)\mT$, where all stability matrices $\bar B_n,\tilde B_{n+1}\T$ are bounded by one with the possible exception of the two end steps, which do not spoil the estimates.
Some smoothness of the problem is required, and with the abbreviations $z=(y\T,p\T,u\T)\T$ and $f(z)=\big(f_i(z)\big)_{i=1}^m:=f(y,u)$, we assume bounds
\begin{align}\label{LipVoraus}
 \|\nabla_z f(z)\|,(1+\|z\|)\|\nabla_{zz} f(z)\|\le\Lambda_2,
 \; \|\langle p,\nabla_{uzz}f(z)\rangle\|\le\Lambda_3,
\end{align}
in a tubular neighborhood of the exact solution $z^\star(t)$.
For the sake of clarity for derivatives including summations in $\R^m$, we introduced the notation $\langle p,\nabla_z f\rangle:=\sum_{i=1}^m(e_i\T p)\nabla_z f_i$, etc.
Since $K$ is a diagonal matrix in all steps now, the strengthened {\em control uniqueness property} \cite{Hager2000} needs to be assumed at the exact solution only as
\begin{align}\label{defntuu}
 \|\langle p^\star(t),\nabla_{uu}f(z^\star(t))\rangle^{-1}\|\le\omega,\ t\in[0,T].
\end{align}
The Peer methods exploit super-convergence, consequently, the global errors depend on more than one derivative of the solution.
These estimates will be abbreviated by a Matlab-type notation
\begin{align}\label{ablnrm}
 \|\phi^{(q:r)}\|_{[n]}:=\max_{q\le k\le r}\|\phi^{(k)}\|_{[n]},\quad 0\le q\le r,\ 0\le n\le N,
\end{align}
for sufficiently smooth functions $\phi(t)$.
Only linear objective functions $\CC$ are considered in order to shorten the proof.
\begin{theorem}\label{TGlobFe}
Let the solution $z^\star=\big((y^\star)\T,(p^\star)\T,(u^\star)\T\big)\T$ of \eqref{KKT_state}--\eqref{KKT_ctr} be smooth with $y^\star,p^\star\in C^{q+1}[0,T]$, let the objective function $\CC$ in \eqref{OCprob_objfunc} be linear and $f$ satisfy the assumptions \eqref{LipVoraus} in some tubular neighborhood of $z^\star$ and \eqref{defntuu} in $z^\star$ itself.
Assume that the Peer triplet satisfies the order conditions \eqref{OBed_vw}, \eqref{OBed_ad} and \eqref{OBedstrt}, \eqref{OBedend} with $q\ge1$.
Let there exist a fixed weight matrix $W\in\R^{s\times s}$ such that $\|W^{-1}A^{-1}B(\sigma)W\|_\infty=1$ and $\|W^{-1}A^{-1}B(\sigma)W-e_1e_1\T\|_\infty\le\tilde\gamma<1$ according to \eqref{decoup} for $\sigma\in[\underline\sigma,\bar\sigma]$ with $0<\underline\sigma<1<\bar\sigma$.
Then, under the condition
\begin{align}\label{Liprestr}
 \zeta_1 T\Lambda_2\max\{1,4\omega\Lambda_2\}\le\frac23,
\end{align}
where $\zeta_1$ is a constant only depending on the Peer triplet, there exists a unique solution $Z=(Y\T,P\T,U\T)\T$ to the system \eqref{KKT_state_peer_init}--\eqref{KKT_ctr_peer} for grids $\{t_n\}$ with stepsize ratios $\sigma_n\in[\underline\sigma,\bar\sigma]$ and  $H=\max_{n=0}^Nh_n$ small enough.
Under the additional condition \eqref{supkvvg} for super-convergence, the error to the exact solution ${\bf z}_n\T=({\bf y}_n\T,{\bf p}_n\T,{\bf u}_n\T)$ is bounded as
\begin{align}\notag
 \max_{n=0}^N\{\|Y_n-{\bf y}_n\|_\infty,&\|P_n-{\bf p}_n\|_\infty,\|U_n-{\bf u}_n\|_\infty\}\\
 &\le\mu\max_{n=0}^N h_n^{q}\max\big\{\|(y^\star)^{(q:q+1)}\|_{[n]},\|(p^\star)^{(q:q+1)}\|_{[n]}\big\}.\label{globFe}
\end{align}
with some constant $\mu$.
The estimate \eqref{globFe} also holds under the weaker condition \eqref{supkvvs} for smooth grids obeying $\sigma_n=1+O(h_n)$, see \eqref{glatt}.
\end{theorem}
We note that the Theorem applies to the Peer triplets \texttt{AP4o33vgi}, \texttt{AP4o33vsi} with $q=s-1=3$ and that the error bound relates the local stepsize to the local size of solution derivatives as a justification for the kind of error control discussed in the next section.\\
\par\noindent
\textbf{Sketch of proof} similar to \cite{LangSchmitt2024}, Theorem~7.2:
Errors are denoted as $\check Z_n=Z_n-{\bf z}_n$.
Subtracting in a usual way the equations \eqref{KKT_state_peer}, \eqref{KKT_adj_peer} defining the standard method, and those for the truncation error, \eqref{locerry}, \eqref{locerrp}, we get for the standard method the error equations
\begin{align}\label{fegly}
 \check Y_n-\bar B_n\check Y_{n-1}&=R_n^Y(\check Z)-\tau_n^Y,\;n=1,\ldots,N,\\\nonumber
 R_n(\check Z)&:=h_n\bar K_n\big(F({\bf z}_n+\check Z_n)-F({\bf z}_n)\big),\\\label{feglp}
 \check P_n-\tilde B_{n+1}\T\check P_{n+1}&=R_n^P(\check Z_n)-\tau_n^P,\;n=0,\ldots,N-1,\\
  R_n^P(\check Z)&:=h_nA_n\mT\big(\langle K({\bf p}_n+\check P_n),
  \nabla_YF({\bf z}_n+\check Z_n)\rangle-\langle K{\bf p}_n,\nabla_YF({\bf z}_n\rangle\big),\nonumber
\end{align}
with $\bar B_n=A_n^{-1}B_n$, $\bar K_n=A_n^{-1}K$, $\tilde B_{n+1}\T=A_n\mT B_{n+1}\T$.
In the equations \eqref{fegly}, \eqref{feglp}, the linear parts with constant coefficients were collected on the left hand sides.
The error equations for the starting step \eqref{KKT_state_peer_init} and the adjoint boundary condition \eqref{KKT_adj_peer_init} are treated in a similar manner.
The constraint \eqref{KKT_ctr_peer} with $N_U=\{0\}$ reads $\kappa_{ii}\langle P_{ni},\nabla_u f(Z_{ni})\rangle=0$, $n=0,\ldots,N$, $i=1,\ldots,s$, and has exactly the same form as \eqref{KKT_ctr}, which means that there is no truncation error here, i.e., $\tau_n^U=0$, and no smoothness assumptions on the control are required.
In order to obtain the same structure for this constraint with constant coefficients on the left hand side as in \eqref{fegly}, \eqref{feglp}, we have to fix derivatives at the exact solution ${\bf z}_n$.
Also $\kappa_{ii}>0$ may be dropped in these conditions.
So for the total global error $\check Z=\big(\check Z_n\big)_{n=0}^N$, we end up with the system
\begin{align}\label{gesfe}
\MM\check Z=-\tau+R(\check Z):=\begin{pmatrix}
 -\tau^Y+R^Y(\check Z)\\
 -\tau^P+R^P(\check Z)\\
 R^U(\check Z)
\end{pmatrix},\quad
\MM=\begin{pmatrix}
 M_{11}\otimes I_m&0&0\\
 0&M_{22}\otimes I_m&0\\
 \MM_{31}&\MM_{32}&\Omega
\end{pmatrix},
\end{align}
where $\MM_{31}=\mbox{blockdiag}_{n,i}\big(\langle {\bf p}_{ni},\nabla_{uy}f({\bf z}_{ni})\rangle\big)$, $\MM_{32}=\mbox{blockdiag}_{n,i}\big(\nabla_{u}F({\bf z}_{ni})\T \big)$, and in the main diagonal, $\Omega=\mbox{blockdiag}_{n,i}\big(\langle{\bf p}_{ni},\nabla_{uu}f({\bf z}_{ni})\rangle\big)$.
Protracted computations in \cite{LangSchmitt2024} showed Lipschitz conditions for the right-hand sides $R_n$ of the form
\begin{align}\label{LipRYP}
 \|R_n^Y(\tilde Z)-R_n^Y(\hat Z)\|,\,\|R_n^P(\tilde Z)-R_n^P(\hat Z)\|
 \le h_n\zeta_1\Lambda_2\|\tilde Z_n-\hat Z_n\|,\\\label{LipRU}
\|R_n^U(\tilde Z)-R_n^U(\hat Z)\|
 \le (\zeta_2\Lambda_2+\zeta_3\Lambda_3)(\|\tilde Z_n\|+\|\hat Z_n\|)\|\tilde Z_n-\hat Z_n\|,
\end{align}
for $\tilde Z_n,\hat Z_n$ in neighborhoods of the origin.
We note that the constants $\zeta_i$ only depend on the coefficients of the Peer method.
The matrices $M_{11},M_{22}$ in \eqref{gesfe} are block bidiagonal matrices with identities $I_s$ in the main diagonal.
According to \eqref{fegly}, $M_{11}$ has $\bar B_n$ in its first subdiagonal blocks while the first superdiagonal blocks in $M_{22}$ are $\tilde B_{n+1}\T$ by \eqref{feglp}.
We recall that both satisfy $\ltnorm\bar B_n\rtnorm=\ltnorm\tilde B_{n+1}\rtnorm=1$ in the respective weighted norms.
Hence, rewriting \eqref{gesfe} in fixed point form $\check Z=\Phi(\check Z)$ with $\Phi(\check Z)=\MM^{-1}(R(\check Z)-\tau)$, we obtain in these weighted norms for the first two parts $\Phi^Y,\Phi^P$ the estimates
\begin{align}\label{LipPhYP}
 \ltnorm\Phi^Y(\tilde Z)-\Phi^Y(\hat Z)\rtnorm,\,\ltnorm \Phi^P(\tilde Z)-\Phi^P(\hat Z)\rtnorm
 \le T\zeta_1\Lambda_2\ltnorm\tilde Z-\hat Z\rtnorm,
\end{align}
where $\ltnorm Z\rtnorm=\max_{n=0}^N\ltnorm Z_n\rtnorm$.
Now, the Lipschitz condition for the last row of \eqref{gesfe} may be rewritten as
\begin{align*}
 \Phi^U(\tilde Z)-\Phi^U(\hat Z)=\Omega^{-1}\Big(-\MM_{31},-\MM_{32},I\Big)
 \begin{pmatrix}
 \Phi^Y(\tilde Z)-\Phi^Y(\hat Z)\\
 \Phi^P(\tilde Z)-\Phi^P(\hat Z)\\
 R^U(\tilde Z)-R^U(\hat Z)
 \end{pmatrix},
\end{align*}
leading with assumption \eqref{defntuu} and \eqref{LipPhYP} to
\begin{align*}
 \ltnorm \Phi^U(\tilde Z)-\Phi^U(\hat Z)\rtnorm\le
 2\omega((\zeta_2\Lambda_2+\zeta_3\Lambda_3)\varepsilon+\zeta_1T\Lambda_2^2\big)
 \ltnorm \tilde Z-\hat Z\rtnorm
\end{align*}
for $\ltnorm\tilde Z\rtnorm,\,\ltnorm\hat Z\rtnorm\le\varepsilon$.
Choosing a small neighborhood of the origin with radius $\hat\varepsilon>0$ such that $(\zeta_2\Lambda_2+\zeta_3\Lambda_3)\hat\varepsilon\le \zeta_1T\Lambda_2^2$, we see that $L_{\Phi}=\zeta_1T\Lambda_2\max\{1,4\omega\Lambda_2\}$ is a Lipschitz constant for $\Phi$ which satisfies $L_\Phi\le2/3$ under assumption \eqref{Liprestr}.
And by standard arguments, the $\hat\varepsilon$-neighborhood is mapped onto itself if $\ltnorm \Phi(0)\rtnorm=\ltnorm\MM^{-1}\tau\rtnorm\le\hat\varepsilon/3$.
Hence, there exists a unique fixed point $\check Z$ with
\begin{align}\label{fixpkt}
\ltnorm\check Z\rtnorm=3\ltnorm\Phi(\check Z)\rtnorm-2\ltnorm\check Z\rtnorm
\le3\ltnorm\Phi(\check Z)-\Phi(0)\rtnorm-2\ltnorm\check Z\rtnorm+3\ltnorm\Phi(0)\rtnorm
\le3\ltnorm\MM^{-1}\tau\rtnorm
\end{align}
under the assumption $3\ltnorm\MM^{-1}\tau\rtnorm\le\hat\varepsilon$.
\par
Applying simple summation of bounds in $\ltnorm\MM^{-1}\tau\rtnorm$ would lose one order of convergence.
This loss can be avoided with the aid of the super-convergence conditions \eqref{supkvvg}, \eqref{supkvvs} along the lines 
of the proof of Lemma 7.3 in \cite{LangSchmitt2024}.
According to \eqref{decoup}, long products $\bar B_n\cdots\bar B_{k-1},\,k<n,$ converge to the rank one matrix $\eins\eins\T A$ as
$\|W^{-1}\bar B_n\cdots \bar B_{k+1}W -e_1e_1\T\|_\infty\le\tilde\gamma^{n-k}$.
Now, due to $\eins\T A\beta_3(\sigma)=e_1\T W^{-1}\beta_3(\sigma)\equiv 0$, the contribution of, e.g., $M_{11}^{-1}\tau^Y$ in \eqref{fixpkt} for time step $n$ is
\begin{align}\notag
 &\|W^{-1}\sum_{k=0}^n\bar B_n\cdots \bar B_{k+1}\tau_k^Y\|_\infty
 \le\|\sum_{k=0}^ne_1e_1\T W^{-1}\tau_k^Y\|_\infty+\sum_{k=0}^n \tilde\gamma^{n-k}\|\tau_k^Y\|_\infty\\\label{Mnulltau}
 &=O\big(\sum_{k=0}^n (h_k^{q+1}\|(y^\star)^{(q+1)}\|_{[n]}+\tilde\gamma^{n-k}h_k^{q}\|(y^\star)^{(q)}\|_{[n]})\big)
 =O\big(\max_{k=0}^n h_k^{q}\|(y^\star)^{(q:q+1)}\|_{[k]}\big).
\end{align}
Hence, for $q\ge1$ and $H$ small enough, the estimate \eqref{Mnulltau} will be smaller than $\hat\varepsilon/3$ and \eqref{fixpkt} yields the assertion.
For the adjoint error $M_{22}^{-1}\tau^P$ of the Peer method \texttt{AP4o33vsi}, where $\beta_3^\dagger$ satisfies \eqref{supkvvs} only, the leading error in $e_1\T(W^\dagger)^{-1}\tau_k^P$ is now $O(h_k^{q+1}\|(p^\star)^{(q:q+1)}\|_{[n]})$ for smooth grids, which is covered by the final bound \eqref{globFe}.
Here, $W^\dagger=(AW)\mT$ is the weight matrix for all matrices $\tilde B_{n+1}\T$.
\qed

\section{A posteriori equi-distribution of global errors}\label{sec:equierr}
In \cite[Chapter 7.1.2]{LangSchmitt2025a}, we have discussed the opportunity to provide an optimized stepsize
sequence a priori using the exact solution $y(t)$ to implement
an error equi-distribution principle introduced
in \cite[Chapter 9.1.1]{AscherMattheijRussell1995}. Here, we are going to exploit a posteriori error estimators
to approximate global errors.
Since the last row of the inverse of a Vandermonde matrix contains the coefficients of the highest difference expression, we set
for our four-stage third-order Peer triplets
\begin{align*}
 v_1\T =&\,6\,e_4\T V_4^{-1},\;
 v_2\T =6\,e_4\T V_4^{-1}\PP_4^{-1},
\end{align*}
and define the estimators as a weighted mean of the polynomial derivatives from two intervals,
\begin{align*}
\varepsilon^Y_0=v_1\T Y_0,\quad
&\varepsilon^Y_n=\delta v_1\T Y_n +
(1-\delta)\sigma_n^{s-1}v_2\T Y_{n-1},\;n=1,\ldots,N,\\[2mm]
\varepsilon^P_N=v_1\T P_N,\quad
&\varepsilon^P_{n-1}=(1-\delta)v_1\T P_{n} +
\delta\sigma_n^{s-1}v_2\T P_{n-1},\;n=1,\ldots,N.
\end{align*}
where $\delta\in [0,1]$ is chosen as a weighting factor. 
In long-standing numerical experiments with Peer methods, we have discovered 
that the use of previous stage values, i.e., $\delta=0$, works quite
reliable for stiff and very stiff problems. For mildly stiff problems, the choice $\delta=1$ often leads to a slightly better
performance. For our examples in Section~\ref{sec:num}, we will present results for $\delta=0$.
\par
Taylor expansions with exact solutions $(y^\star(t),p^\star(t))$ at
off-step points $t=t_{ni}$ yield $\varepsilon^Y_n=h_n^3(y^\star)^{(3)}(t_n)+O(h_n^4)$ and
$\varepsilon^P_n=h_n^3(p^\star)^{(3)}(t_n)+O(h_n^4)$. Supposing an asymptotic behaviour of the global errors,
\[ \| y(t_n)-Y(t_n) \|=O\left( h_n^3\|(y^\star)^{(3)}(t_n)\|\right),\;
\| p(t_n)-P(t_n) \|=O\left( h_n^3\|(p^\star)^{(3)}(t_n)\|\right) \]
with numerical approximations $Y(t_n)=e_1\T V_4^{-1}Y_n$ and $P(t_n)=e_1\T V_4^{-1} P_n$, the estimators defined above can be
used to construct a new time mesh with nearly equally distributed global errors. In order to simplify
equi-distribution of these errors, often linearized versions such as $h_n\|(y^\star)^{(3)}(t_n)\|^{1/3}$ and
$h_n\|(p^\star)^{(3)}(t_n)\|^{1/3}$ are considered, see e.g. \cite{HuangRussell2011}. Defining the extended vector of
approximate third derivatives $d_n:=((\varepsilon^Y_n)\T,(\varepsilon^P_n)\T)\T /h_n^3$ and a piecewise constant
mesh density function $\psi(t) := \|d_n\|^{1/3},\,t\in [t_{n},t_{n+1})$, on the old grid it holds that
\begin{align*}
\int_{t_n}^{t_{n+1}} \psi(t)\,dt = h_n\psi(t_n) = h_n\, \|d_n\|^{1/3} \approx C
\left( \|y(t_n)-Y(t_n)\|^2 + \|p(t_n)-P(t_n)\|^2 \right)^{1/6}.
\end{align*}
The equi-distribution of the global errors over a new mesh $0 < t_1' < . . . < t_{N+1}' = T$ requires that
\begin{align}\label{equdist}
\int_{t'_n}^{t'_{n+1}} \psi(t)\,dt =
 \frac{1}{N+1} \int_0^T\psi(t)\,dt = const,\ n=0,\ldots,N.
\end{align}
Since the global errors for the state and adjoint variables can differ by several orders of magnitude,
we have to apply an appropriate weighting procedure to make sure that both errors are equally considered in
the mesh design. Note that improving the approximate adjoints is a key to enhance the quality of the
computed control. This motivates the definition of the following practical, component-wise weighted measures:
\begin{align*}
\theta^Y_n:=&\,err_{3,n}\max_{i=1,\ldots,m}
\frac{|\varepsilon^Y_{ni}|}{atol_Y+rtol_Y\,\hat{Y}_{ni}},\;n=0,\ldots,N,\\
\theta^P_n:=&\,err_{3,n}^\dagger\max_{i=1,\ldots,m}
\frac{|\varepsilon^P_{ni}|}{atol_P+rtol_P\,\hat{P}_{ni}},\;n=0,\ldots,N,
\end{align*}
with the intermediate values
\begin{align*}
\hat{Y}_0=&\, |Y(t_0)|,\;\hat{Y}_n=\delta |Y(t_n)| + (1-\delta) |Y(t_{n-1})|,\;n=1,\ldots,N,\\[2mm]
\hat{P}_N=&\, |P(t_N)|,\;\hat{P}_{n-1}=\delta |P(t_{n-1})| + (1-\delta) |P(t_n)|,\;n=1,\ldots,N,
\end{align*}
and the individual error constants $err_{3,n}$, $err_{3,n}^\dagger$ from Table~\ref{TPT} and Table~\ref{TPB}.
Selecting appropriate values for the scaling parameters $atol_Y$, $atol_P$, $rtol_Y$, and $rtol_P$ is crucial 
and should reflect the nature of the problem.
Defining the vector of the balanced weighted error estimators
\begin{align*}
\theta_n:=\left( (\theta_n^Y)\T,\omega(\theta_n^P)\T\right),\;
\omega=\|\theta_n^Y\|_\infty/\|\theta_n^P\|_\infty,
\end{align*}
we eventually define the mesh density function by $\psi(t):=\|\theta_n/h_n^3\|^{1/3}$, $t\in [t_n,t_{n+1})$.
\par
Then, an efficient way to construct a new mesh $0=t'_0<t'_1<\ldots<t'_{N+1}=T$ which distributes
a given mesh density function $\psi(t)>0$ evenly according to (\ref{equdist})
is the computation of a smooth mesh function $x(\xi),\,\xi\in[0,1]$, which defines the grid through
\begin{align}\label{gitxxi}
 t_n'=x(\xi_n),\ \xi_n=\frac{n}{N+1},\ n=0,\ldots,N+1.
\end{align}
According to \cite[Chapter~2.2.2]{HuangRussell2011} such a mesh function may be computed as the
solution of the nonlinear boundary value problem
\begin{align}\label{gitrwp}
 \left(\psi(x) x_\xi\right)_\xi=0,\quad x(0)=0,\ x(1)=T,
\end{align}
where the lower index $\xi$ denotes derivation.
Since the full order of some Peer triplets requires that the used grids are smooth satisfying
$\sigma_n=1+O(h_n)$, we will discuss now how the expression
\begin{align}\label{defeta}
 \eta_n:=\frac{\sigma_n-1}{h_n}=\frac{h_n-h_{n-1}}{h_{n-1}h_n}
\end{align}
relates to the grid functions $x(\xi),\psi(t)$.
\begin{lemma}\label{Leta}
Let the grid $\{t_n'\}$ be defined according to \eqref{gitxxi} with some smooth mesh function $x(\xi)$,
$h_n':=t_{n+1}'-t_n'$, and $\sigma_n':=h_{n}'/h_{n-1}'$.
Then,
\begin{align}\label{etax}
\eta_n'=\frac{\sigma_n'-1}{h_n'}
 =\frac{x_{\xi\xi}(\xi_n)}{(x_\xi(\xi_n))^2}+O\left( \frac1{(N+1)^{2}} \right).
\end{align}
If the mesh function $x(\xi)$ is the solution of the boundary value problem \eqref{gitrwp} with some smooth,
positive density function $\psi$, then also
\begin{align}\label{etapsi}
\eta_n'
 =-\frac{\psi'(t_n')}{\psi(t_n')}+O\left( \frac1{(N+1)^{2}}\right)
 =-\big(\ln\psi\big)'\big|_{t_n'}+O\left( \frac1{(N+1)^{2}}\right).
\end{align}
\end{lemma}
{\bf Proof.} For ease of writing, we abbreviate $\delta:=(N+1)^{-1}$.\\
a) Using Taylor expansion at $\xi_n$ in \eqref{defeta}, we get
\begin{align*}
h_n'-h_{n-1}'=&\,x(\xi_{n-1})-2x(\xi_n)+x(\xi_{n+1})=\delta^2x_{\xi\xi}(\xi_n)+\frac{\delta^4}{12}x^{(4)},\\
h_{n-1}'h_n'=&\,\delta^2\big(x_\xi(\xi_n)\big)^2+\frac{\delta^2}{12}(4x_\xi(\xi_n)x^{(3)}-3x_{\xi\xi}(\xi_n)^2),
\end{align*}
where the highest derivatives are evaluated at some intermediate point.
This shows \eqref{etax}.
\\{} b)
The differential equation in \eqref{gitrwp} may also be written as $\psi x_{\xi\xi}+\psi'x_\xi^2=0$, which shows that $x_{\xi\xi}/x_\xi^2=-\psi'/\psi$.\qed\\
\par
For the solution of (\ref{gitrwp}), we use MMPDELab by W.~Huang \cite{Huang2019} and start with a uniform grid.
Linear finite elements with $N$ inner nodes $\xi_n$, $n=1,\ldots,N$, and a pseudo-timestepping scheme are applied.
Additionally, we apply the smoothing of the mesh density function $\psi(t)$ in MMPDELab by controlling the ratios
$\psi'(t_n')/\psi(t_n')$ from (\ref{etapsi}) in order to guarantee
$\sigma_n'=1+\eta_n'h_n'$ with $|\eta_n'|\le 15$ on the new mesh $t_n'=x(\xi_n)$. Such meshes yield
smoothly varying stepsize sequences and work quite well for our variable stepsize Peer methods.
Good starting values for the new control vector $U'$ at $\{t_{ni}'\}$ can be computed using
piecewise cubic polynomial interpolation of $U$. A pseudocode for the overall algorithm with the key parameters 
is given in Table~\ref{ALG}.
\par
In what follows, we will apply this approach to two problems: (i) boundary control for
the 1D heat equation with known exact discrete solutions \cite{LangSchmitt2023} and (ii) optimal
control of cytotoxic therapies on prostate cancer growth in 2D \cite{ColliEtAl2021}. All calculations
have been done with Matlab-Version R2019a on a Latitude 7280 with an i5-7300U Intel processor at 2.7 GHz.
\begin{table}
\centering
\begin{tabular}{p{0.5cm}p{10.5cm}}
\hline\noalign{\smallskip}
\multicolumn{2}{l}{Algorithm: Adaptive Time Grids}\\
\noalign{\smallskip}\hline\noalign{\smallskip}
1. & Solve optimal control problem on uniform grid $\{t_n\}$ for control $U$.\\[1mm]
2. & Given weight $\delta\in [0,1]$, compute estimators $\varepsilon^Y_i$, $\varepsilon^P_i$, $i=0,\ldots,N$, 
     of the global errors.\\[1mm]
3. & Given scaling parameters $atol_Y$, $atol_P$, $rtol_Y$, and $rtol_P$, build mesh density function $\psi(t)$.\\[1mm]
4. & Given smoothness parameter $\eta>0$, solve mesh PDE with smoothing of $\psi(t)$ such that
     $\sigma'_n=1+\eta'_nh'_n$ with $|\eta'_n|\le\eta$ for the new grid $\{t'_n\}$.\\[1mm]
5. & Update grid from $\{t_n\}$ to $\{t'_n\}$.\\[1mm]
6. & Calculate initial guess for $U'$ at $\{t_{ni}'\}$ by interpolating $U$.\\[1mm]
7. & Solve optimal control problem on grid $\{t'_n\}$.\\
\noalign{\smallskip}\hline\noalign{\smallskip}
\end{tabular}
\parbox{14cm}{
\caption{Algorithm for the adaptive time-grid construction by 
a posteriori equi-distribution of global errors.\label{ALG}}
}
\end{table}
\section{Numerical results}\label{sec:num}
\subsection{Boundary control for the 1D Heat equation}
The first problem, originally introduced in \cite{LangSchmitt2023}, was specifically constructed to
yield exact analytical solutions for an optimal boundary control problem governed
by a one-dimensional discrete heat equation. The objective function evaluates both
the deviation of the final state from a desired target and the cost of the control.
Since there are no spatial discretization errors, the numerical convergence orders
of time integration methods can be determined with high precision, without the need
for a numerical reference solution.
\par
The optimal control problem reads as follows:
\begin{align*}
\mbox{minimize } &\,\frac12\|y(1)-\hat{y}\|^2_2 +
\frac{1}{2}\int_{0}^{1}u(t)^2\,dt\\[1mm]
\mbox{subject to } &\,y'(t) = Ay(t) +  \gamma e_m u(t),\quad t\in (0,1],\\
&\, y(0) = \eins_m,
\end{align*}
with
\begin{align*}
 A =&\frac1{(\triangle x)^2}\begin{pmatrix}
   -1&1\\
   1&-2&1\\
   &&\ddots&\ddots&\ddots\\
   &&&1&-2&1\\
   &&&&1&-3
 \end{pmatrix},
\end{align*}
state vector $y(t)\in\R^m$, $\triangle x=1/m$, and
$\gamma=2/(\triangle x)^2$. We set $m=250$. The components $y_i(t)$
approximate the solution of the continuous 1D heat equation $Y(x,t)$ over the spatial domain
$[0,1]$ at the discrete points $x_i=(i-0.5)\triangle x$, $i=1,\ldots,m$. The corresponding
boundary conditions are $\partial_xY(0,t)=0$ and $Y(1,t)=u(t)$. The matrix $A\in\R^{m\times m}$
results from standard central finite differences. Its eigenvalues $\lambda_k$ and
corresponding normalized orthogonal eigenvectors $v^{[k]}$ are given by
\begin{align*}
\lambda_k =&\;-4m^2\sin^2\left( \frac{\omega_k}{2m}\right),\;
\omega_k=\left( k-\frac12 \right)\pi,\\[2mm]
v^{[k]}_i=&\;\nu_k\cos\left( \omega_k\frac{2i-1}{2m}\right),\;
\nu_k = \frac{2}{\sqrt{2m+\sin(2\omega_k)/\sin(\omega_k/m)}},\;\; i,k=1,\ldots,m.
\end{align*}
Introducing an additional
component $y_{m+1}(t)$ and adding the equations $y'_{m+1}(t)=u(t)^2$, $y_{m+1}(0)=0$,
the objective function can be transformed to the Mayer form
\[ \CC(y(1)):=\frac12\,\left( \sum_{i=1}^{m}(y_i(1)-\hat{y}_i)^2+y_{m+1}(1)\right)\]
with the extended vector $\bar y(1)=(y_1(1),\ldots,y_m(1),y_{m+1}(1))\T$.
\par
Following the test case in \cite{LangSchmitt2023}, we define the target vector
$\hat{y}$ through
\begin{align*}
\hat{y}(t) =&\,y^\star(T)-\delta\left( v^{[1]} + v^{[2]}\right).
\end{align*}
with $\delta=1/75$ and the exact solution $y^\star(T)=\sum_{i=1,\ldots,m}\eta_k(T)v^{[k]}$.
The coefficients $\eta_k(T)$ are given by
\begin{align*}
\eta_k(T) =&\;e^{\lambda_kT}\eta_k(0)-\gamma^2\delta Tv_m^{[k]}\sum_{l=1}^{2}
v_m^{[l]}\varphi_1((\lambda_k+\lambda_l)T)
\end{align*}
where $\eta_k(0)=(y^\star(0))\T v^{[k]}$ and $\varphi_1(z):=(e^z-1)/z$.
The exact control and adjoint are
\begin{align*}
u^\star(t) = -\gamma p_m(t),\quad
p^\star(t) =&\; \delta \left( e^{\lambda_1(T-t)}v^{[1]} +  e^{\lambda_2(T-t)}v^{[2]}\right).
\end{align*}
We will compare the numerical errors for $y(T)$, $p(0)$ and $u(t)$.
An approximation $p_h(0)$ for the Peer method is obtained from $p_h(0)=e_1\T V_4^{-1}P_0$.
Note that, compared to \cite{LangSchmitt2023}, we have changed the sign of the
adjoint variables, i.e., $p\mapsto -p$, to fit into our setting.

\begin{figure}[t!]
\centering
\includegraphics[width=7cm]{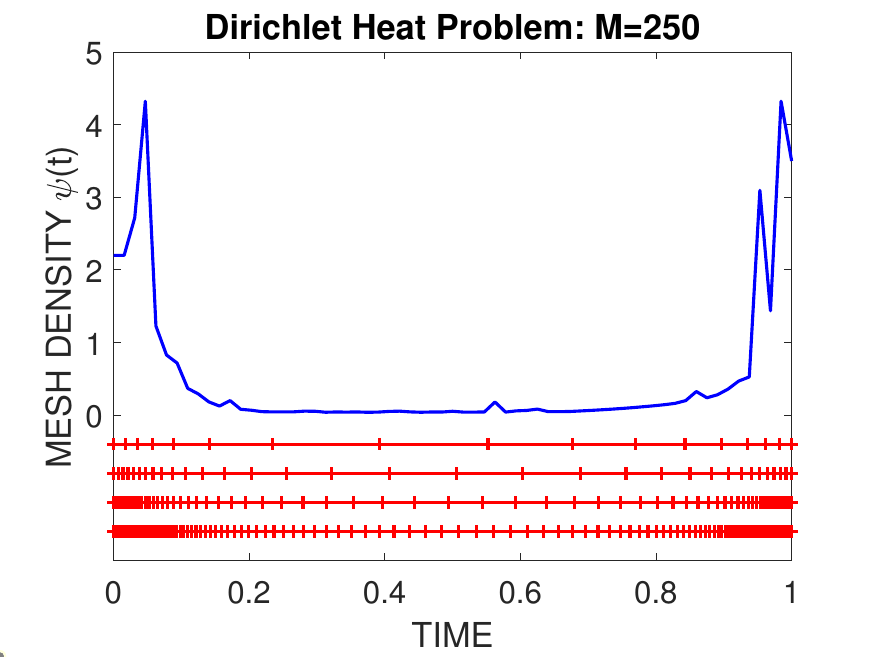}
\includegraphics[width=7cm]{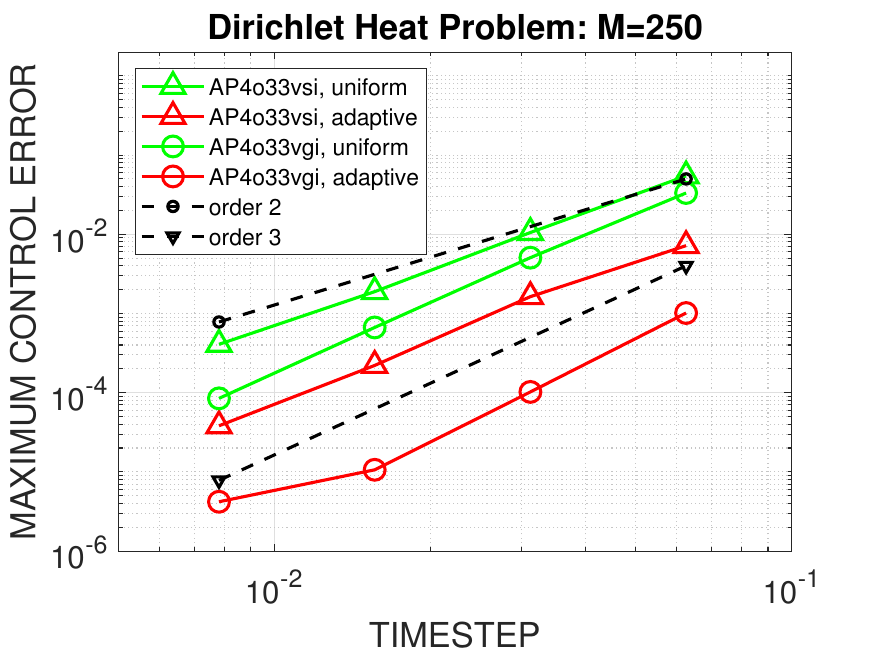}\\[5mm]
\includegraphics[width=7cm]{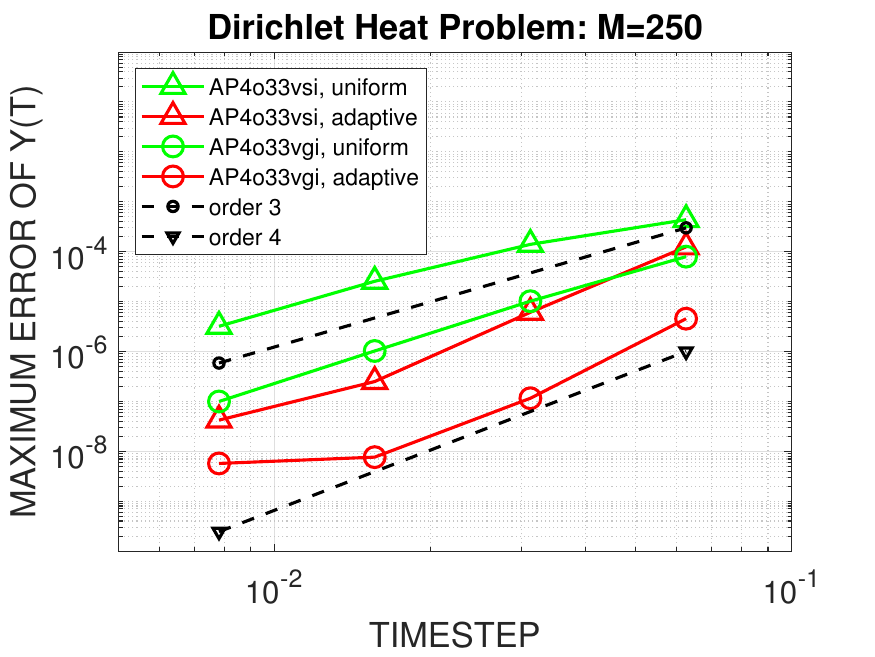}
\includegraphics[width=7cm]{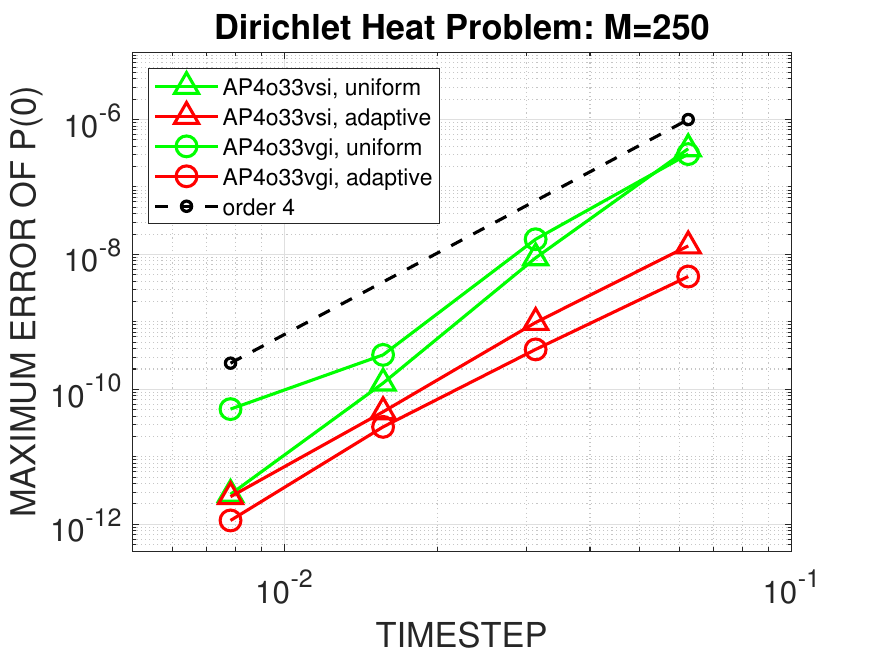}
\parbox{13.5cm}{
\caption{Dirichlet heat problem with $m=250$ spatial points. Exemplary mesh density
function $\psi(t)$
for \texttt{AP4o33vgi}, $N=63$ and adapted time grids for $N=15,31,63,123$ (top left). Convergence of
the maximal control errors $\|U_{ni}-u(t_{ni})\|_\infty$ (top right), state errors
$\|y(T)$$-$$y_h(T)\|_\infty$ (bottom left), and adjoint errors
$\|p(0)$$-$$p_h(0)\|_\infty$ (bottom right) for uniform and adaptive time grids.
\label{fig:hproblem}}
}
\end{figure}
We performed calculations for $N+1=16,32,64,128$ with the results shown in Figure 1.
The default interior-point algorithm in \texttt{fmincon} was used, with a zero control
vector as the initial guess for uniform time grids. 
We set $\delta=0$, $atol_Y=atol_P=10^{-8}$, 
and $rtol_Y=rtol_P=1$ to calculate the mesh density function. The self-adjoint method
\texttt{AP4o33vgi} performs remarkably well, achieving average convergence orders of
$3.2$ for the state and $4.2$ for the adjoint variables. This results in third-order
convergence for the controls, even on coarse time grids.
For variable-stepsizes, observed orders of convergence exceeding the expected ones may be 
explained by the fact that on fine grids a disproportionate amount of points may be 
placed in the critical parts of the problem domain.
The method \texttt{AP4o33vsi}, which offers a larger angle for $A(\alpha)$-stability
but has fewer symmetry properties, attains an asymptotic convergence order of $3$ for
the state and an unexpectedly high average of $5.7$ for the adjoint. This yields an
average convergence order of $2.4$ for the control.
\par
The mesh density function $\psi(t)$ reflects the expected behaviour for the linear
heat equation. Larger, well-balanced global errors are observed at both
ends of the time interval (see Figure~\ref{fig:hproblem}, top left). Solving
\eqref{gitrwp} using MMPDELab with smoothing yields improved variable time grids with
significant refinements near the boundaries. The resulting meshes satisfy
$\sigma_n=1+\eta_nh_n$ with $|\eta_n|\le 15$.
More precisely, we get
\begin{align*}
\min_n{\sigma_n} =
\left\{
 \begin{array}{lll}
 (0.69,0.76,0.87,0.96) &\;\quad \text{for } \texttt{AP4o33vgi},\\[2mm]
 (0.65,0.77,0.88,0.94) &\;\quad \text{for } \texttt{AP4o33vsi},\
 \end{array}
\right.\\[2mm]
\max_n{\sigma_n} =
\left\{
 \begin{array}{lll}
 (1.75,1.29,1.15,1.03) &\;\quad \text{for } \texttt{AP4o33vgi},\\[2mm]
 (1.33,1.31,1.14,1.06) &\;\quad \text{for } \texttt{AP4o33vsi}.\
 \end{array}
\right.
\end{align*}
The improvement in control accuracy achieved by both methods is significant,
approximately by a factor of $10$ for \texttt{AP4o33vsi} and nearly $50$ for
\texttt{AP4o33vgi}, except at $N+1=64$, where the improvement is somewhat lower.
This clearly demonstrates the potential of using adaptive time grids and time
integrators that maintain their order for such grids. In all cases, the post-processed
control variables $U^\ddagger_{ni}$ defined in \eqref{ctrpp} do not show a 
visible improvement and therefore do not contribute to a higher convergence order. 
As a result, we omit the details.
\par
The triangular approximation at the boundaries, i.e., 
using $\tilde A_i$ instead of $A_i$,
$i=0,N$, requires $10-15$ iteration steps of the form \eqref{Nwtit} to meet a relative 
tolerance of $10^{-14}$ in the maximum norm for all four smaller subsystems
which are solved by the sparse direct solver of {\sc Matlab}. For practical tolerances as 
$10^{-6}$, only $5-7$ iteration steps are necessary. This results in approximately
$25$ linear direct solves of size $m=250$ for each boundary step instead of 
one direct solve for the fully coupled (linear) systems written out as $(A_i\otimes I_m-h_iK\otimes{\bf J}_i)Y_i=b_i$ and
$(A_i\T\otimes I_m-h_iK\otimes {\bf J}_i\T) P_i={b}_i^\dagger$, $i=0,N$, of size
$4*m=1,000$ with full matrices $A_i\in\R^{4\times 4}$.
\par
\begin{remark}
It is a known fact that one-step methods such as Rosenbrock, Runge–Kutta,
and extrapolation methods suffer from order reduction when applied to
stiff ODEs or semi-discrete nonlinear parabolic PDEs 
\cite{LubichOstermann1995,OstermannRoche1992}. In the latter case, the temporal 
order of convergence is mainly influenced by the spatial regularity of the solution, 
which usually depends on the boundary conditions. Time-dependent boundary conditions 
constitute the most severe case. Thus, classical symplectic integrators,
such as the A-stable 
Runge-Kutta-Gauss methods, lose their advantage of high-order accuracy when 
applied to boundary control problems. For example, the two-stage fourth-order 
Runge-Kutta-Gauss method
shows a reduction to first order for the approximation of the control. 
The interested reader is referred to
\cite{LangSchmitt2023} for a detailed discussion of this phenomenon. 
\end{remark}
\subsection{2D Prostate cancer growth}
The second problem, taken from \cite{ColliEtAl2021}, addresses the optimal control of
cytotoxic and antiangiogenic therapies in the treatment of prostate cancer (PCa).
The mathematical model in a three-dimensional space-time cylinder $Q_T:=\Omega\times (0,T)$
with boundary $\Sigma_T:=\partial\Omega\times (0,T)$ reads as follows:
\begin{align}
\partial_t\phi - \lambda\triangle\phi &=\,
-F'(\phi) + (m(\sigma)-U)h'(\phi) &\text{in } Q_T,\\[2mm]
\partial_t\sigma - \eta\triangle\sigma &=\,
-\gamma_h\sigma - (\gamma_c-\gamma_h)\sigma\phi+S_h(1-\phi)+(S_c-S)\phi &\text{in }Q_T,\\[2mm]
\partial_tp - D\triangle p &=\,-\gamma_pp+\alpha_h+(\alpha_c-\alpha_h)\phi &\text{in }Q_T,
\end{align}
with boundary and initial conditions
\begin{align}
\phi=0,\;\partial_n\sigma=\partial_np=0 \quad\text{in }\Sigma_T,\\[2mm]
\phi(0,x)=\phi_0(x),\;\sigma(0,x)=\,\sigma_0(x),\;p(0,x)=p_0(x) \quad\text{in }\Omega.
\end{align}
Here, $\phi$ denotes the phase field variable that identifies the spatial location of the tumour,
$\sigma$ is the concentration of vital nutrients, and $p$ is the tissue PSA concentration
(prostate specific antigen). The controls $U$ and $S$ describe the cytotoxic and antiangiogenic
treatment effects, respectively. In the phase field equation it holds
\[ F(\phi) = M\phi^2(1-\phi)^2,\; h(\phi)=M\phi^2(3-2\phi), \]
with mobility $M>0$. The function $m(\sigma)$ is defined as
\[ m(\sigma) = m_{ref} \left( \frac{\rho+A}{2} + \frac{\rho-A}{\pi}
\arctan \left( \frac{\sigma-\sigma_l}{\sigma_r}\right)\right),\]
where $m_{ref}>0$. $A$ and $\rho$ are constants that determine the rates of cell proliferation
and cell death. We set $\rho=K_\rho/\bar{K}_\rho$ and $A=-K_A/\bar{K}_A$.
The constants $\sigma_l$ and $\sigma_r$ are threshold and reference values for the nutrients,
respectively. A characterization of all parameters and their values which are used in our computation
are given in Appendix~C.
\par
The value for the clinically used serum PSA is determined from
\[ P_s(t) = \int_\Omega\,p(t,x)\,dx. \]
Following \cite{ColliEtAl2021}, we define the two-dimensional simulation domain $\Omega$ as
a square with edge length $l_d=3000\,\si{\micro\meter}$.
The initial tumour is selected ellipsoidally in the centre of the domain:
\[ \phi_0(x) = 0.5 - 0.5 \tanh\left( 10\left(
\sqrt{\frac{(x_1-l_d/2)^2}{a_1^2} + \frac{(x_2-l_d/2)^2}{a_2^2}}-1\right)\right) \]
with $a_1=150\,\si{\micro\meter}$ and $a_2=200\,\si{\micro\meter}$. The initial values for the
nutrients and PSA are estimated as
\[\sigma_0(x) = c_\sigma^0+c_\sigma^1\phi_0(x),\quad p_0(x) = c_p^0+c_p^1\phi_0(x)\]
with $c_\sigma^0=1\,\si{\gram}/\si{\litre}$,
$c_\sigma^1=-0.8\,\si{\gram}/\si{\litre}$, $c_p^0=0.0625\,\si{\nano\gram}/\si{\milli\litre}/\mathrm{cm}^3$ and
$c_p^1=0.7975\,\si{\nano\gram}/\si{\milli\litre}/\mathrm{cm}^3$.
\par
Initially the tumor growth is untreated for $60$ days, i.e., $U(t)=0$ and $S(t)=0$.
The results from \cite{ColliEtAl2021} show that strong cytotoxic chemotherapy suffices to
optimally control the volume and serum PSA of the tumor. So we also set $S(t)=0$ for
later times. The values of $\phi$, $\sigma$, and $p$ obtained at the end of the
untreated growth phase are used as initial conditions for the optimal control problem with
drug delivery. The time is reset to $t=0$ and $T=21\,\text{days}$ is used as a common duration of
a therapy cycle for PCa.
\par
The authors in \cite{ColliEtAl2021} apply $256$ isogeometric elements per side of the computational domain
and a second-order time integrator with uniform stepsize $h=0.1\,\text{days}$ to derive an
optimal single cycle with smoothly decreasing $U_{d1}(t)$, starting with an initial guess
$U_0(t)=m_{ref}\beta_cd_ce^{-t/\tau_c}$, where $\tau_c=5\,\text{days}$,
$\beta_c=1.59\,10^{-2}\,\si{\meter}^2/\si{\milli\gram}$, $d_c=75\,\si{\milli\gram}/\si{\meter}^2$
\cite[Table~2]{ColliEtAl2021}. The maximum value is set to $U_{max}=0.12/\text{day}$, which
corresponds to a $100\,\si{\milli\gram}/\si{\meter}^2$ drug dose of docetaxel. Running a nonlinear
least-square fit of protocols in a post-processing treatment, a new 3-dose drug design of the form
\[ U_{d3}(t) = \sum_{i=1}^{3} m_{ref}\beta_c d_{c,i} e^{-\frac{t-t_{c,i}}{\tau_c}}H(t-t_{c,i}) \]
has been proposed with
$d_c=(58.49,9.20,5.03)\,\si{\milli\gram}/\si{\meter}^2$,
$t_c=(2.85,7.90,9.16)\,\text{days}$ \cite[Table~3]{ColliEtAl2021}, and $H$ being the Heaviside function.
The objective functional is given by
\begin{align*}
J(\phi,p,U) =&\, k_1\int_{Q_T} \phi(t,x)^2\,dxdt + k_2\int_\Omega \phi(T,x)^2\,dx \\[1mm]
&\;+ k_3\int_0^T \left( P_s(t) - P_\Omega \right)^2\,dt
+ k_4\int_{Q_T}\,(U-U_d)^2\,dxdt,
\end{align*}
where $P_\Omega=\alpha_h|\Omega|/\gamma_c$, $U_d=0$ and the weights were set to $k_1=k_2=k_3=1$, $k_4=0.5$.
Note that $p(x,t)-\alpha_h/\gamma_c\ge 0$ due to $p_0(x)\ge\alpha_h/\gamma_p$ and $\phi(x,t)\ge 0$
in the model equations for $p$.
\par
The use of a uniform stepsize of $h=0.1\,\text{days}$, which requests the numerical solution
of $3650$ systems of nonlinear equations with hundred thousands of
spatial degrees of freedoms for each forward and backward run, and the two-step nature of the calculation
of optimal controls is a serious limitation of the proposed algorithm. This becomes even more evident when
drug protocols are to be designed for 3D PCa growth scenarios and longer simulation times are taken into
account, reflecting conditions more closely aligned with experimental and clinical settings. A few
numerical experiments with our third-order Peer triplets already showed that around $100$ variable
time steps are sufficient to reproduce the results above with acceptable accuracy.
We now want to document this in detail.
\par
Let $Z_m:=\{x_i=(i-0.5)\triangle x,\;i=1,\ldots,m\}$ be a set of discrete points lying in $[0,l_d]$. Then
standard central differences of the diffusion
terms and the homogeneous Neumann boundary conditions yield a second-order finite difference approximation
$y(t)\T:=(\Phi(t)\T,\Sigma(t)\T,P(t)\T)\in \R^{3m^2}$ of the PCa model equations on the two-dimensional grid 
$Z\times Z$ with long stacked space grid vectors $\Phi(t)\in\R^{m^2}$, etc.
We use $m=256$, that means $\triangle x=11.72\,\si{\micro\meter}$, which results in using $6$ discrete points to approximate the
interface defined by the distance between $\phi=0.05$ and $\phi=0.95$. Applying row-wise numbering of
the spatial degrees of freedom, the
spatial integrals of the objective functional are discretized by the second-order trapezoidal rule with a row vector of weights
\[ W_{tp} = \frac{(\triangle x)^2}{4}(w_1,w_2,\ldots,w_2,w_1)\in\R^{m^2}, \]
where
\[ w_1=(1,2,\ldots,2,1)\in\R^m,\quad
  w_2=(2,4,\ldots,4,2)\in\R^m. \]
Then the discrete objective functional reads
\begin{align*}
J_{\triangle x}(\Phi,P,U) =&\, k_1 \int_0^T W_{tp} \Phi(t)^2\, dt + k_2 W_{tp} \Phi(T)^2 \\[1mm]
&\, + k_3 \int_0^T \left( W_{tp} \left( P(t) - \frac{\alpha_h}{\gamma_p}\eins \right)\right)^2\,dt
+ k_4 \int_0^T W_{tp}\eins\, (U(t)-U_d(t))^2\,dt
\end{align*}
with $\eins\in\R^{m^2}$. A transformation to Mayer form is done by introducing a further
component $y_{3m^2+1}$ that satisfies
\begin{align*}
y'_{3m^2+1}(t) &=\, k_1 W_{tp} \Phi(t)^2 + k_3 \left( W_{tp} \left( P(t) -
\frac{\alpha_h}{\gamma_p}\eins \right)\right)^2 + k_4 W_{tp}\eins\, (U(t)-U_d(t))^2,\\
y_{3m^2+1}(0) &=\, 0.
\end{align*}
Eventually, we get the reduced discrete objective function
\[ C(y(T)) = k_2 W_{tp} \left( \left( y(T)\right)_{1}^{m^2}\right)^2 + y_{3m^2+1}(T).\]
The admissible set of controls is $U_{ad}:=\{U(t)\in\R:0\le U(t)\le U_{max}\}$.
We solved the optimal 
control problem with \texttt{AP4o33vgi} using $N=83$ time steps
of variable size. The algorithm \textit{trust-region-reflective} was applied as recommended
large-scale solver for \texttt{fmincon}, where $U_0(t)$ was taken as initial guess. 
The scaling parameters for the calculation of the
mesh density function were set to $atol_Y=10^{-8}$, $atol_P=10^2$, and $rtol_Y=rtol_P=1$.
Once again, we used $\delta=0$.
\begin{figure}[t!]
\centering
\includegraphics[width=7cm,height=5cm]{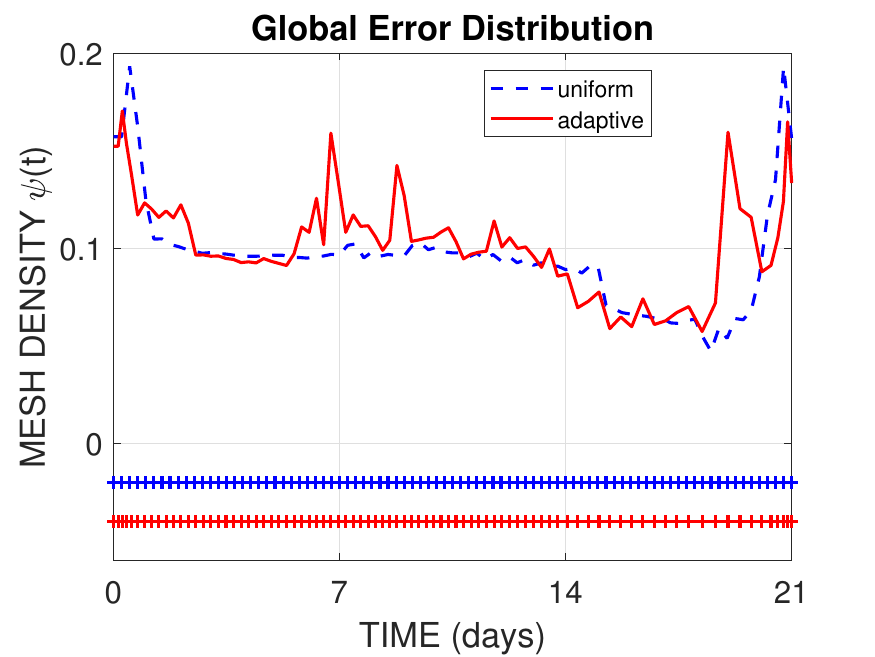}
\includegraphics[width=7cm,height=5cm]{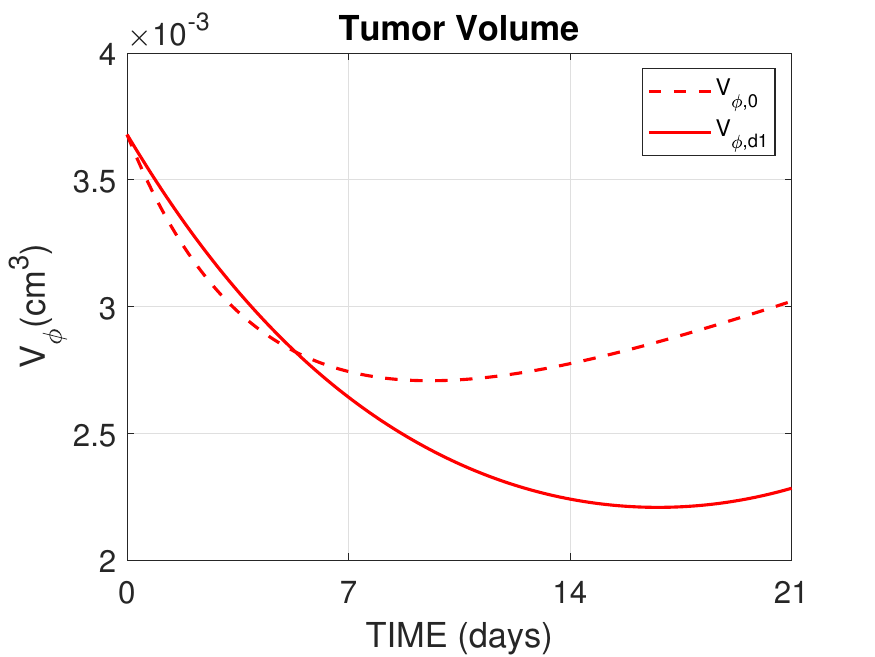}\\[5mm]
\includegraphics[width=7cm,height=5cm]{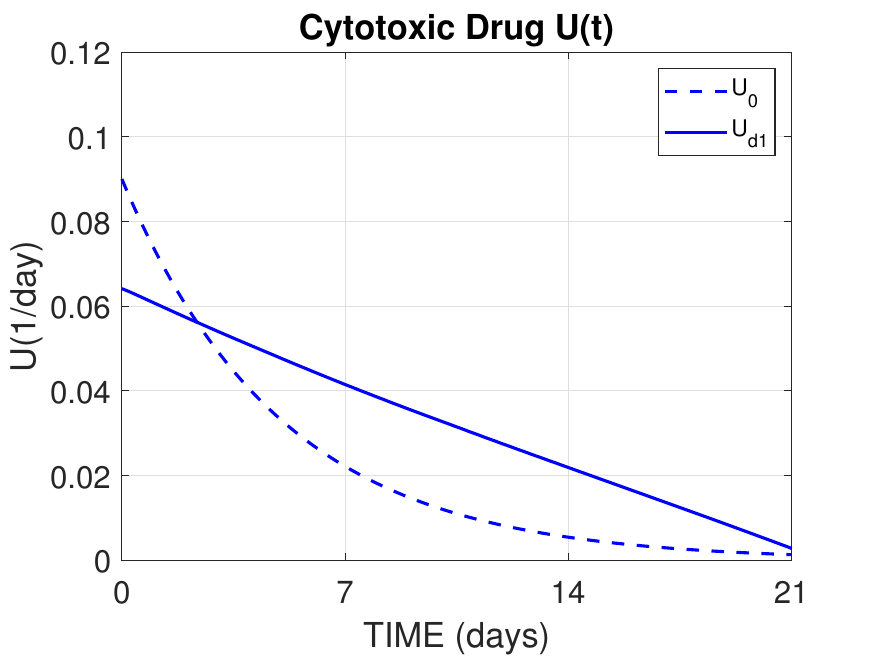}
\includegraphics[width=7cm,height=5cm]{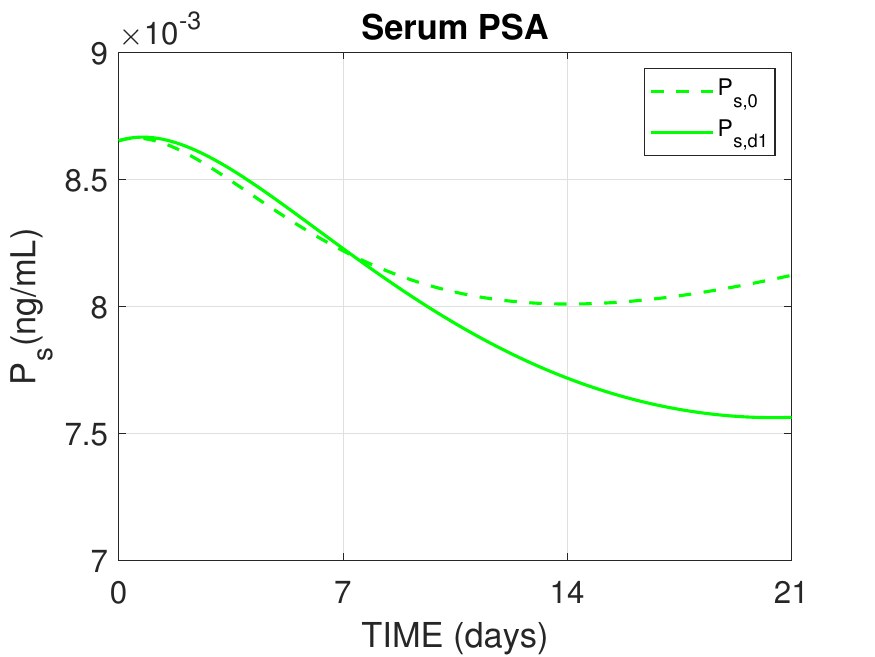}
\parbox{13.5cm}{
\caption{\texttt{AP4o33vgi} for PCa problem: single drug cycle. The target optimal 
cytotoxic 1-dose drug $U_{d1}(t)$ 
(bottom left) yields a significant smaller tumor volume $V_{\phi,d1}$ (top right) and 
serum PSA $P_{s,d1}$ (bottom right), compared to $V_{\phi,0}$ and $P_{s,0}$ obtained 
for the standard docetaxel protocol $U_0(t)$. Refining the time grid at both ends of the 
interval using the equi-distribution principle reduces the mesh density function there (top left)
and leads to a reduction of the estimated maximum global 
error by nearly 40\%. The measures of the adaptive mesh are $\sigma_n\in [0.75,1.26]$ 
and $|\eta_n|\le 2$.
\label{fig:pca_run1}}
}
\end{figure}
Two different drug protocols were considered. In the first test, we reproduced the optimal 
single cycle with smoothly decreasing $U_{d1}(t)$ from \cite{ColliEtAl2021} setting $U_d=0$ 
in the objective function. Choosing $k_1=k_2=k_3=1$ and $k_4=6$, we achieved a good agreement with the
results discussed in \cite[Figure~4]{ColliEtAl2021} with much fewer time steps.
The results are presented in Figure~\ref{fig:pca_run1}. The \texttt{fmincon} solver converges in just 4 iterations, 
achieving 4-digit accuracy in the objective function within 27 minutes.
\begin{figure}[t!]
\centering
\includegraphics[width=7cm,height=5cm]{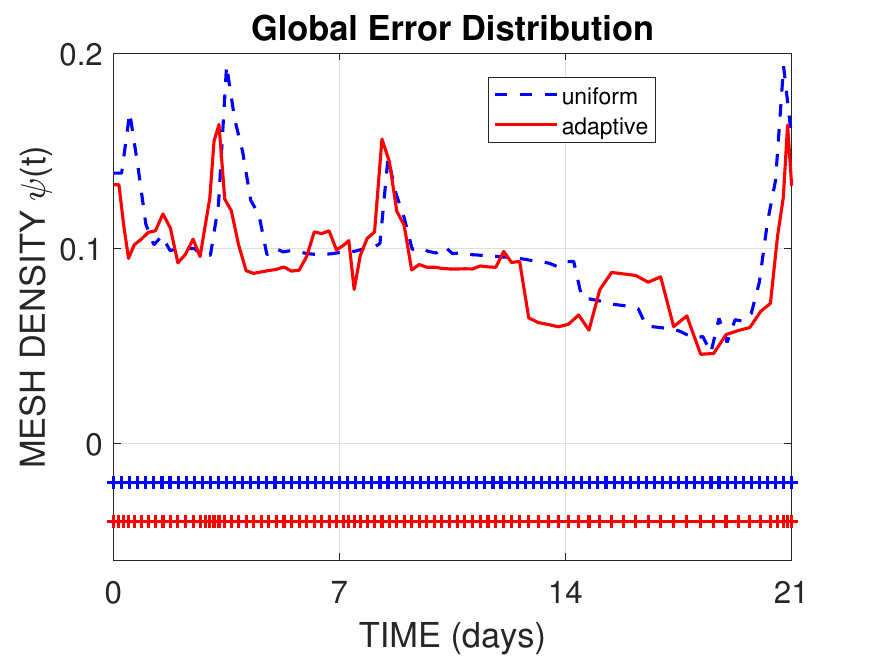}
\includegraphics[width=7cm,height=5cm]{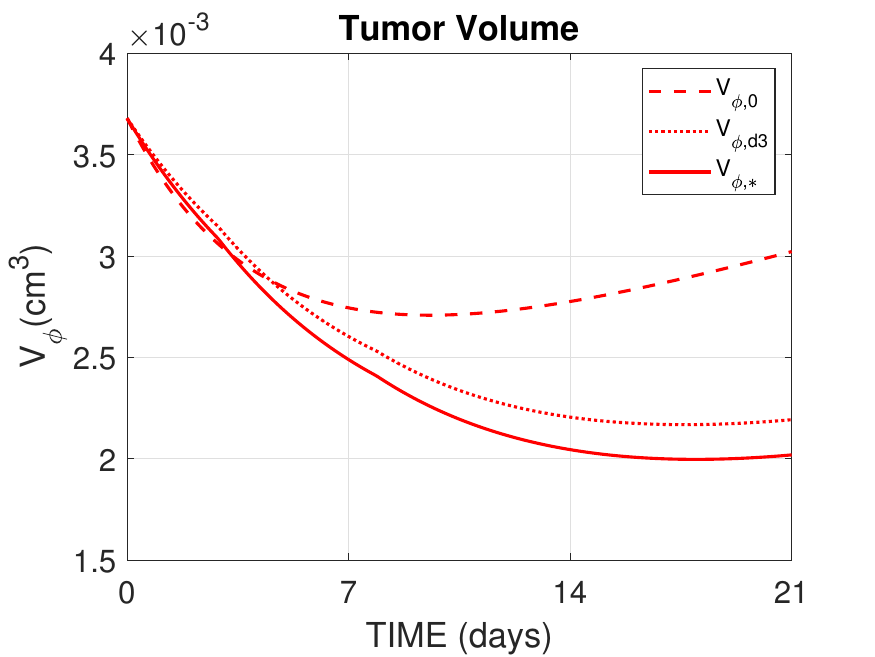}\\[5mm]
\includegraphics[width=7cm,height=5cm]{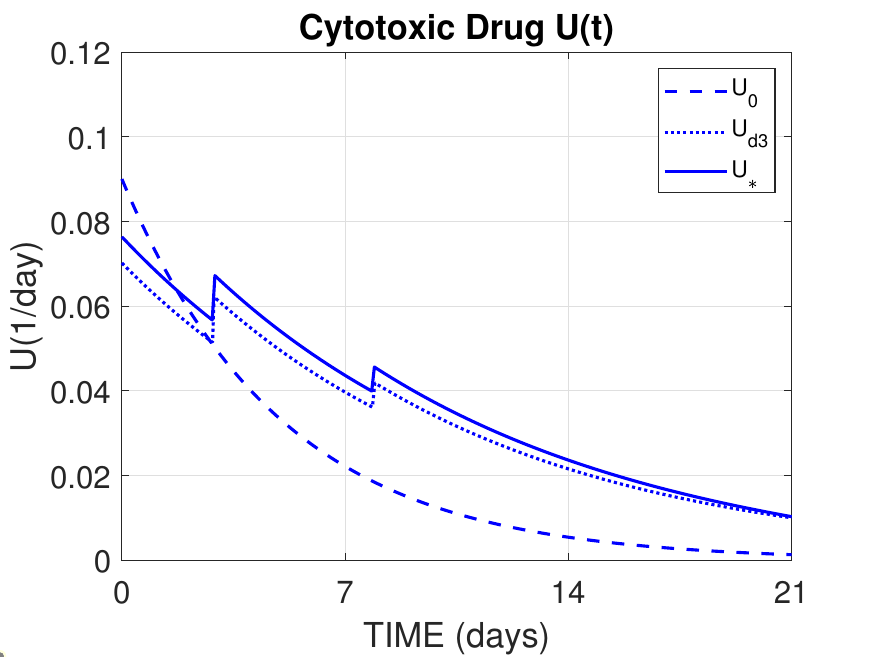}
\includegraphics[width=7cm,height=5cm]{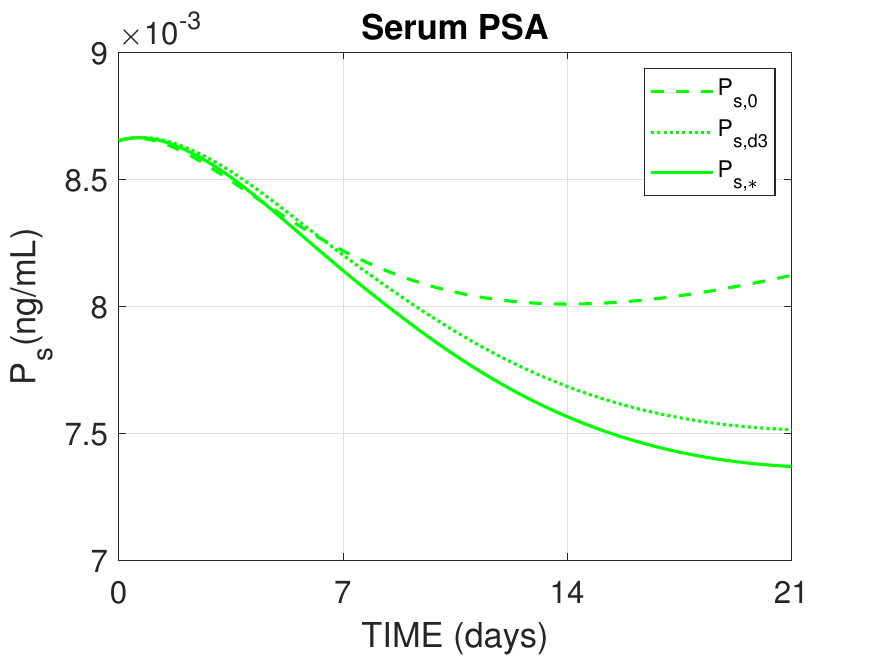}
\parbox{13.5cm}{
\caption{\texttt{AP4o33vgi} for PCa problem: 3-dose drug cycle. The numerical target 
optimal cytotoxic drug $U_\star(t)$
(bottom left) yields a reduced tumor volume $V_{\phi,\star}$ (top right) and 
serum PSA $P_{s,\star}$ (bottom right), compared to $V_{\phi,d3}$ and $P_{s,d3}$ obtained 
for the 3-dose docetaxel protocol $U_{d3}(t)$. Refining the time grid at both ends of the 
interval and at the drug delivery time points $t_c$ using the equi-distribution principle 
reduces the mesh density function in critical regions (top left) and leads to a reduction of
the estimated maximum global error by 38\%. The measures of the adaptive 
mesh are $\sigma_n\in [0.75,1.20]$ and $|\eta_n|\le 2$.
\label{fig:pca_run2}}
}
\end{figure}
\par%
The second numerical study was concerned with the newly designed 3-dose drug protocol
$U_{d3}(t)$ from above. Here, we set $U_d=U_{d3}$ in the objective function and 
investigated the evolution of the tumor volume and the serum PSA for different choices
of the penalty parameter $k_4$. In Figure~\ref{fig:pca_run2}, we summarize the results
for $k_4=60$.
A slight increase in the docetaxel doses improves the therapeutic result.
Our algorithm is directly applicable, and the discontinuities in drug administration are 
resolved with sufficient accuracy. The quality of the computed gradients is remarkably
high. So \texttt{fmincon} again achieves convergence in 4 iterations only, with 6-digit 
accuracy in the objective function. In both tests, the post-processed control variables 
$U^\ddagger_{ni}$ do not improve the already high, third-order accuracy of the approximate 
controls.
\par
By using triangular approximations at the boundaries, i.e., replacing $A_i$ by $\tilde A_i$,
$i=0,N$, and solving linear systems of the form \eqref{Nwtit}, it is still 
possible to apply a sparse direct solver to linear systems of size $196,609$. Typical 
iteration counts are around 12 for a relative tolerance of $10^{-3}$, while approximately 
$6$ Newton iterations are typically required for the four nonlinear subsystems of
the standard Peer steps with $i=1,\ldots,N-1$. This means that the strong coupling of the
stage solutions at the boundaries resulting in nonlinear systems of size $786,436$ can be 
resolved by solving a sequence of smaller systems of size $196,609$ at
the expense of only doubling the computational effort in the boundary steps.
\par\noindent
\section{Conclusion}\label{sec:con}
We modified two of the third-order Peer two-step triplets from \cite{LangSchmitt2025a} and integrated
automatic grid construction in efficient gradient-based solution algorithms for the fully coupled
optimal control problem with unknown control $u(t)$. Positivity of the quadrature 
weights in all matrices $K$ is essential. To avoid additional memory requirements in large nonlinear systems, 
special triangular iterations with small contraction factors are applied in the initial and final steps. 
The pulcherrima triplet \texttt{AP4o33vgi} possesses LSRK-form \eqref{LSRK} for state and 
adjoint equations, yielding super-convergence order three for arbitrary stepsize ratios. An additional
grid-smoothness condition $|\sigma_n-1|\le \eta h_n$ is necessary to ensure third-order adjoint 
super-convergence for the Peer triplet \texttt{AP4o33vsi}. Extensive numerical experiments showed that 
$\eta\le 15$ is an appropriate choice for the design of smooth adaptive time grids. This moderate 
restriction comes with an significantly increased stability angle of $83.74^o$, compared to $61.59^o$
of the pulcherrima triplet.
\par
Both schemes, \texttt{AP4o33vgi} and \texttt{AP4o33vsi}, perform very well on a benchmark problem involving 
boundary heat control — an application where one-step methods typically suffer from order reduction, as 
demonstrated in \cite{LangSchmitt2023}. Noteworthy, the symmetry properties of the pulcherrima scheme 
\texttt{AP4o33vgi} yield increased orders between three and four for all components $(y,u,p)$, even on
coarse grids. 
\par
As a second example, we applied our solution strategy to a real-world optimal control problem in medical 
treatment planning for prostate cancer, as discussed in \cite{ColliEtAl2021}. In contrast to the main 
discretization parameters used in that work, our approach requires only 83 time steps of variable size — 
compared to 3650 for the generalized $\alpha$-scheme — and just a few conjugate gradient steps using the 
large-scale optimization solver \textit{trust-region-reflective} from {\sc Matlab}. This replaces the simpler 
\textit{steepest-descent} gradient algorithm with step-size control and a fixed maximum of 100 iterations, 
while still achieving high-quality solutions. It also reflects the high-order consistency of the gradient
provided by the self-adjoint \texttt{AP4o33vgi}. A further advantage is the opportunity to
approximate optimal controls with steep gradients, allowing the direct study of alternative medical treatments 
close to real multi-dose drug protocols. 
\par
In both applications, the post-processed control $U^\ddagger$ introduced in \eqref{ctrpp} offers no discernible
benefit. We infer that for Peer methods, which perform close to their
theoretical order or even better, the approximation quality of $U$ is nearly optimal and postprocessing
is not advisable in general.
\par
Summarizing, we believe that the two newly developed Peer two-step methods have the potential to significantly enhance the efficiency of solving large-scale optimal control problems. This includes applications such as designing effective drug protocols for 3D prostate cancer growth scenarios in clinical settings.

\vspace{0.5cm}
\par
\noindent {\bf Acknowledgements.}
The first author is supported by the Deutsche Forschungsgemeinschaft
(German Research Foundation) within the collaborative research center
TRR154 {\em ``Mathemati\-cal modeling, simulation and optimisation using
the example of gas networks''} (Project-ID 239904186, TRR154/3-2022, TP B01).

\begin{appendix}
\section{Coefficients of \texttt{AP4o33vgi} and \texttt{AP4o33vsi}}\label{SApendx}
The essential data are the node vector $\cc\T=(c_1,c_2,c_3,c_4)$, the coefficients
$(A_0,K)$, $(A,K)$, $(A_N,K)$ of the starting, the standard and the end step and the free
coefficients of the sparse matrix
\begin{align*}
 \hat B(\sigma)=\begin{pmatrix}
  1&1&1&1\\
  0&0&0&\hat b_{24}(\sigma)\\
  0&0&0&\hat b_{34}(\sigma)\\
  \hat a_{41}&\hat b_{42}(\sigma)&\hat b_{43}(\sigma)&\hat b_{44}(\sigma)
 \end{pmatrix},
\end{align*}
depending on the stepsize ratio $\sigma$. All other coefficients may be computed by
\begin{align*}
 a=A_0\eins,\ w=A_N\T\eins,\quad B(\sigma)=V_4\mT\hat B(\sigma)V_4,
\end{align*}
where $V_4=(\eins,\cc,\cc^2,\cc^3)$ is the Vandermonde matrix for the nodes in $\cc$.
Since the matrices $A_0,A_N$ from the boundary steps are not in lower triangular form,
we provide additional lower triangular approximations $\tilde A_0,\tilde A_N$ to be used
in fast Gauss-Seidel type block iterations \eqref{Nwtit} in order to avoid the solution
of large coupled systems of fourfold size.
Since the elements of  $\tilde A_0,\tilde A_N$ coincide with those of  $A_0,A_N$ in all
subdiagonals, only the diagonals of $\tilde A_0,\tilde A_N$ are given.

\subsection*{A1: Coefficients of \texttt{AP4o33vgi}}

\[ \cc\T=\left(0,\frac13,\frac23,1 \right),\quad
  K=\diag\left(\frac18,\frac38,\frac38,\frac18\right),\]
\[
A=\begin{pmatrix}
  1&0&0&0\\[1.5mm]
  -\frac94&\frac94&0&0\\[1.5mm]
  \frac94&-\frac92&\frac94&0\\[1.5mm]
  -1&\frac94&-\frac94&1
 \end{pmatrix},\quad
 \hat B(\sigma)=\begin{pmatrix}
  1&1&1&1\\[1mm]
  0&0&0&\frac1{36\sigma}\\[1mm]
  0&0&0&0\\[1mm]
  0&\frac{\sigma}{36}&\frac{\sigma}{18}&
  \frac1{804}(132\sigma+\frac{65}\sigma-149)
 \end{pmatrix}
 \]
 \[
A_0=\begin{pmatrix}
 \frac{47161}{23112}&   \frac{945}{1712}& \frac{9}{856}& -\frac{113}{1712}\\[1.5mm]
 -\frac{41383}{7704}&  \frac{1017}{1712}& -\frac{27}{856}& \frac{339}{1712}\\[1.5mm]
  \frac{41383}{7704}& -\frac{4869}{1712}& \frac{1953}{856}& -\frac{339}{1712}\\[1.5mm]
 -\frac{47161}{23112}& \frac{2907}{1712}& -\frac{1935}{856}& \frac{1825}{1712}
 \end{pmatrix},\;
 A_N=\begin{pmatrix}
 \frac{1825}{1712}& -\frac{339}{1712}& \frac{339}{1712}& -\frac{113}{1712}\\[1.5mm]
 -\frac{1935}{856}& \frac{1953}{856}& -\frac{27}{856}& \frac{9}{856}\\[1.5mm]
 \frac{2907}{1712}& -\frac{4869}{1712}& \frac{1017}{1712}& \frac{945}{1712}\\[1.5mm]
 -\frac{47161}{23112}& \frac{41383}{7704}& -\frac{41383}{7704}& \frac{47161}{23112}
\end{pmatrix},
\]
\[ \left(\tilde a_{ii}^{(0)}\right)=\left(\frac{154}{75},\frac{69}{40},\frac{219}{94},\frac{67}{63}\right),\quad
 \left(\tilde a_{ii}^{(N)}\right)
 =\left(\frac{67}{63},\frac{219}{94},\frac{69}{40},\frac{154}{75}\right).
\]
\subsection*{A2: Coefficients of \texttt{AP4o33vsi}}
\begin{align*}
 \cc\T=&\left(\frac{144997}{389708},\frac{73}{748},\frac{77297572}{117896267},1\right),\\[1mm]
 K=&\diag\left(0.2089552772313791,0.2461266069992848,
  0.4259606950456414,0.1189574207236947\right),
\end{align*}
\begin{align*}
A_0=&\begin{pmatrix}
 1.26852968140859992& -2.79702966259295784&  0.0151774841161155076& 0\\
0.254440961986028910&  1.58797813851094452& -0.00536671649536513773& 0\\
-3.75232398970999177&  2.14140637287657549&  2.46031830832026582& 0\\
 2.22935334631536294& -0.932354848794562167& -2.47012907594101619& 1
\end{pmatrix},\\
A=&\begin{pmatrix}
 0.7588470158140062&0&0&0\\
 0.4346633458753195& 0.5989561692950702&0&0\\
 -3.295204661275873&-0.3671669165116753& 2.473930545531403&0\\
  2.101694299586548&-0.2317892527833949&-2.473930545531403& 1
\end{pmatrix},\\
A_N=&\begin{pmatrix}
 0.721680741868241430& 0.0131418918926231641& 0.0333333333333333333& -0.00930895128019174555\\
 0.123032993110224916&  0.709147801969229717&  0.279492058866634697& -0.078053338775699573\\
 -1.03159221459763137&  -1.16757403034966595&  0.443763401719389714&  0.566961810971761768\\
  5.56340552222272135&  -1.45584078718664692&  -5.57863709363081650&  1.86704685986649197
\end{pmatrix},
\end{align*}
\begin{align*}
 \big(\tilde a_{ii}^{(0)}\big)=&\big(1.58950617283950617,1.66216216216216216, 2.47,1\big),\\
 \big(\tilde a_{ii}^{(N)}\big)
 =&\big(0.725,0.681818181818181818,2,1.91525423728813559\big),
\end{align*}
\begin{align*}
\hat b_{24}=&0.02321239244678227/\sigma,
\;\hat b_{34}=0,\\
\hat a_{41}=&0.1010743874247749,\;
\hat b_{42}=\hat a_{41}+0.003586671392069201\,\sigma,\\
 \hat b_{43}=&\hat a_{41}+0.007173342784138403\,\sigma-0.002465255918355442\,\sigma^2,\\
 \hat b_{44}=&0.0078782707622298066+0.1683589306029579\,\sigma-0.1125\,\sigma^2+0.025\,\sigma^3,
\end{align*}
\section{List of frequently used symbols}
\begin{tabular}{|l|p{11cm}|}\hline
$A,B(\sigma),K$& coefficients of the standard method, $A$ triangular, $K\succ0$ diagonal\\
$A_0,a,A_N,w$, $\tilde A_0,\tilde A_N$& coefficients of the two boundary methods and triangular approximations for iterative solution\\
$\bar B_n=A_n^{-1}B_n,\,\tilde B_{n+1}\T =A_n\mT B_{n+1}\T$&corresponding zero-stability matrices\\
$Y_n=\big(Y_{ni}\big)_{i=1}^s$, $U_n,P_n$&stacked numerical stages of state, control, and Lagrange multipliers\\[1mm]
$t_n,h_n,\sigma_n\in[\underline\sigma,\bar\sigma]$&gridpoints, stepsizes, and stepsize ratios $\sigma_n=h_n/h_{n-1}$\\[1mm]
$\cc,V_q,\PP_q=\exp(\tilde E_q)$&node vector, $s\times q$-Vandermonde matrix, $q\times q$-Pascal matrix\\[1mm]
$S_{n,q}=\diag_{i=1}^q(\sigma_n^{i-1})$&stretching matrix for Taylor coefficients\\[1mm]
$\QQ_{q,q}=V_q\T B(\sigma)V_q\PP_q^{-1}$&plays central role,  order conditions imply $\QQ_{q,q}=e_1e_1\T$, $q<s$\\[1mm]
$\hat A_n=V_s\T A_nV_s$, $\hat K,\hat B(\sigma)$&congruent matrices having restricted form related to $\QQ_{q,q}$\\[1mm]
$W\in\R^{s\times s}$&fixed weight matrix which block-diagonalizes all $\bar B(\sigma)=A^{-1}B(\sigma)$\\[1mm]
$\tau_n^Y,\tau_n^P$&local errors of state and adjoint variables, leading terms $\beta_{q,n}(\sigma),\beta_{q,n}^\dagger(\sigma)$\\[1mm]
$err_{q,n}$, $err_{q,n}^\dagger$&leading error constants for Peer method ($=\|\beta_{q,n}(1)\|_\infty$), and adjoint\\[1mm]
$\|\phi^{(q:r)}\|_{[n]}$&local supremum norm on 3 adjacent subintervals of  $\phi^{(k)}$, $q\le k\le r$\\[1mm]
$\varepsilon_n^Y,\varepsilon_n^P$, $v_1\T,v_2\T$&a-posteriori error estimates and its coefficients\\[1mm]
$\psi(x)$&mesh density function for grid optimization\\\hline
\end{tabular}
{
\section{Model parameters for prostate cancer growth}
\begin{table}[h!]
\centerline{\begin{tabular}{lccc}
\hline
Parameter & Notation & Value \\\hline
\textbf{Tumor dynamics} && \\
Diffusivity of the tumor phase field & $\lambda$ & $640\,\si{\micro\meter}^2/\text{day}$ \\
Tumor mobility & $M$ &  $2.5/\text{day}$ \\
Net proliferation scaling factor &  $m_{ref}$ & $7.55\,10^{-2}/\text{day}$ \\
Scaling reference for proliferation rate & $\bar{K}_p$ & $1.50\,10^{-2}/\text{day}$ \\
Proliferation rate & $K_p$ & $1.50\,10^{-2}/\text{day}$ \\
Scaling reference for apoptosis rate & $\bar{K}_A$ & $2.10\,10^{-2}/\text{day}$ \\
Apoptosis rate & $K_A$ &  $1.37\,10^{-2}/\text{day}$ \\[1mm]
\textbf{Nutrient dynamics} && \\
Nutrient diffusivity & $\eta$ & $6.4\,10^4\,\si{\micro\meter}^2/\text{day}$ \\
Nutrient supply in healthy tissue & $S_h$ & $2\,\si{\gram}/\si{\litre}/\text{day}$ \\
Nutrient supply in tumor tissue & $S_c$ & $2.75\,\si{\gram}/\si{\litre}/\text{day}$ \\
Nutrient uptake by healthy tissue & $\gamma_h$ & $2\,\si{\gram}/\si{\litre}/\text{day}$ \\
Nutrient uptake by tumor tissue & $\gamma_c$ & $17\,\si{\gram}/\si{\litre}/\text{day}$ \\
Hypoxic-viable threshold for nutrients & $\sigma_l$ & $0.4\,\si{\gram}/\si{\litre}$ \\
Scaling reference for nutrients& $\sigma_r$ & $6.67\,10^{-2}\,\si{\gram}/\si{\litre}$ \\[1mm]
\textbf{Tissue PSA dynamics} && \\
Tissue PSA diffusivity & $D$ & $640\,\si{\micro\meter}^2/\text{day}$ \\
Healthy tissue PSA production rate & $\alpha_h$ & $1.712\,10^{-2}\,
\si{\nano\gram}/\si{\milli\litre}/\mathrm{cm}^3/\text{day}$ \\
Tumoral tissue PSA production rate & $\alpha_c$ & $\alpha_c=15\alpha_h$ \\
Tissue PSA natural decay rate & $\gamma_p$ & $ 0.274/\text{day}$ \\\hline
\end{tabular}}
\label{PCAparam}
\end{table}
}
\end{appendix}

\bibliographystyle{plain}
\bibliography{bibpeeropt}

@ARTICLE{LangSchmitt2022a,
  AUTHOR =       {J.~Lang and B.A.~Schmitt},
  TITLE =        {Discrete adjoint implicit peer methods in optimal control},
  JOURNAL =      {J. Comput. Appl. Math.},
  YEAR =         {2022},
  volume =       {416},
  pages =        {114596},
}

@ARTICLE{LangSchmitt2022b,
  AUTHOR =       {J. Lang and B.A.~Schmitt},
  TITLE =        {Implicit {A}-stable peer triplets for {ODE} constrained optimal control problems},
  JOURNAL =      {Algorithms},
  YEAR =         {2022},
  volume =       {15},
  pages =        {310},
}

@ARTICLE{LangSchmitt2023,
  AUTHOR =       {J. Lang and B.A.~Schmitt},
  TITLE =        {Exact discrete solutions of boundary control problems for the 1{D} heat equation},
  JOURNAL =      {J. Optim. Theory Appl.},
  YEAR =         {2023},
  volume =       {196},
  pages =        {1106--1118},
  DOI =          {10.1007/s10957-022-02154-4},
}

@ARTICLE{LangSchmitt2024,
  AUTHOR =       {J. Lang and B.A.~Schmitt},
  TITLE =        {Implicit peer triplets in gradient-based solution algorithms for
                  {ODE} constrained optimal Ccontrol},
  JOURNAL =      {J. Optim. Theory Appl.},
  YEAR =         {2024},
  volume =       {203},
  pages =        {985–1026},
  DOI =          {10.1007/s10957-024-02541-z},
}

@ARTICLE{LangSchmitt2025a,
  AUTHOR =       {J.~Lang and B.A.~Schmitt},
  TITLE =        {Variable-Stepsize Implicit Peer Triplets in {ODE} Constrained Optimal Control},
  JOURNAL =      {J. Comput. Appl. Math.},
  YEAR =         {2025},
  volume =       {460},
  pages =        {116417},
}

@BOOK{Troutman1996,
  author =       {J.L.~Troutman},
  title =        {Variational Calculus and Optimal Control},
  publisher =    {Springer, New York},
  year =         {1996},
}

@ARTICLE{Hager2000,
  AUTHOR =       {W.W.~Hager},
  TITLE =        {{R}unge-{K}utta methods in optimal control
                  and the transformed adjoint system},
  JOURNAL =      {Numer. Math.},
  YEAR =         {2000},
  volume =       {87},
  pages =        {247--282},
}

@ARTICLE{HagerRostamian1987,
  AUTHOR =       {W.W.~Hager and R.~Rostamian},
  TITLE =        {Optimal coatings, bang-bang controls, and
                  gradient techniques},
  JOURNAL =      {Optimal Control Applications and Methods},
  YEAR =         {1987},
  volume =       {8},
  pages =        {1--20},
}

@ARTICLE{BeigelMommerWirschingBock2014,
  AUTHOR =       {D.~Beigel and M.S.~Mommer and L.~Wirsching and H.G.~Bock},
  TITLE =        {Approximation of weak adjoints by reverse automatic
                  differentiation of {BDF} methods},
  JOURNAL =      {Numer. Math.},
  YEAR =         {2014},
  volume =       {126},
  pages =        {383--412},
}

@ARTICLE{LangVerwer2013,
  AUTHOR =       {J.~Lang and J.G.~Verwer},
  TITLE =        {{W}-Methods in Optimal Control},
  JOURNAL =      {Numer. Math.},
  YEAR =         {2013},
  volume =       {124},
  pages =        {337--360},
}

@ARTICLE{BonnansLaurentVarin2006,
  AUTHOR =       {F.J.~Bonnans and J.~Laurent-Varin},
  TITLE =        {Computation of order conditions for
                  symplectic partitioned {R}unge-{K}utta
                  schemes with application to optimal control},
  JOURNAL =      {Numer. Math.},
  YEAR =         {2006},
  volume =       {103},
  pages =        {1--10},
}

@ARTICLE{HertyPareschiSteffensen2013,
  AUTHOR =       {M.~Herty and L.~Pareschi and S.~Steffensen},
  TITLE =        {Implicit-Explicit {R}unge-{K}utta schemes for
                  numerical discretization of optimal control problems},
  JOURNAL =      {SIAM J. Numer. Anal.},
  YEAR =         {2013},
  volume =       {51},
  pages =        {1875--1899},
}

@ARTICLE{AlbiHertyPareschi2019,
  AUTHOR =       {G.~Albi and M.~Herty and L.~Pareschi},
  TITLE =        {Linear multistep methods for optimal control problems
                  and applications to hyperbolic relaxation systems},
  JOURNAL =      {Applied Mathematics and Computation},
  YEAR =         {2019},
  volume =       {354},
  pages =        {460--477},
}

@ARTICLE{SanzSerna2016,
  AUTHOR =       {J.M.~Sanz-Serna},
  TITLE =        {Symplectic {R}unge–{K}utta schemes for adjoint equations,
                  automatic differentiation, optimal control, and more},
  JOURNAL =      {SIAM Review},
  YEAR =         {2016},
  volume =       {58},
  pages =        {3--33},
}

@ARTICLE{AlmuslimaniVilmart2021,
  AUTHOR =       {I.~Almuslimani and G.~Vilmart},
  TITLE =        {Explicit stabilized integrators for stiff optimal
                  control problems},
  JOURNAL =      {SIAM J. Sci. Comput.},
  YEAR =         {2021},
  volume =       {43},
  pages =        {A721--A743},
}

@ARTICLE{LiuFrank2021,
  AUTHOR =       {X.~Liu and J.~Frank},
  TITLE =        {Symplectic {R}unge-{K}utta discretization of a regularized
                  forward-backward sweep iteration for optimal control problems},
  JOURNAL =      {J. Comput. Appl. Math.},
  YEAR =         {2021},
  volume =       {383},
  pages =        {113133},
}

@ARTICLE{OstermannRoche1992,
  AUTHOR =       {A.~Ostermann and M.~Roche},
  TITLE =        {{R}unge-{K}utta methods for partial differential equations
                  and fractional orders of convergence},
  JOURNAL =      {Math. Comp.},
  YEAR =         {1992},
  volume =       {59},
  pages =        {403--420},
}

@ARTICLE{LubichOstermann1995,
  AUTHOR =       {Ch.~Lubich and A.~Ostermann},
  TITLE =        {{R}unge-{K}utta approximation of quasi-linear parabolic equations},
  JOURNAL =      {Math. Comp.},
  YEAR =         {1995},
  volume =       {64},
  pages =        {601--627},
}

@BOOK{HairerWannerLubich2006,
  author =       {E.~Hairer and G.~Wanner and Ch.~Lubich},
  title =        {Geometric Numerical Integration, Structure-Preserving Algorithms
                  for Ordinary Differential Equations},
  series =       {Springer Series in Computational Mathematic},
  volume =       {31},
  publisher =    {Springer, Heidelberg, Berlin},
  year =         {1970},
}

@BOOK{AscherMattheijRussell1995,
  author =       {U.M.~Ascher and R.M.M.~Mattheij and R.D.~Russell},
  title =        {{N}umerical {S}olution of {B}oundary {V}alue {P}roblems
                  for {O}rdinary {D}ifferential {E}quations},
  publisher =    {Society for Industrial and Applied Mathematics},
  year =         {1995},
}

@BOOK{HuangRussell2011,
  author =       {W.~Huang and R.D.~Russell},
  title =        {{A}daptive {M}oving {M}esh {M}ethods},
  publisher =    {Springer New York, Dordrecht, Heidelberg, London},
  year =         {2011},
}

@MANUAL{Huang2019,
  title =        {{MMPDE}lab, https://github.com/weizhanghuang/{MMPDE}lab},
  author =       {W.~Huang},
  year =         {2019},
}

@ARTICLE{Spellucci1998,
  AUTHOR =       {P.~Spellucci},
  TITLE =        {An {SQP} method for general nonlinear programs using only equality constrained subproblems},
  JOURNAL =      {Math. Program.},
  YEAR =         {1998},
  volume =       {82},
  pages =        {413--448},
}

@ARTICLE{ByrdGilbertNocedal2000,
  AUTHOR =       {R.H.~Byrd and J.C.~Gilbert and J.~Nocedal},
  TITLE =        {A trust region method based on interior point techniques
                  for nonlinear programming},
  JOURNAL =      {Math. Program.},
  YEAR =         {2000},
  volume =       {89},
  pages =        {149--185},
}

@ARTICLE{ColemanLi1994,
  AUTHOR =       {T.F.~Coleman and Y.~Li},
  TITLE =        {On the convergence of reflective {N}ewton methods
                  for large-scale nonlinear minimization subject to bounds},
  JOURNAL =      {Math. Program.},
  YEAR =         {1994},
  volume =       {67},
  pages =        {189--224},
}

@ARTICLE{ColliEtAl2021,
  AUTHOR =       {P.~Colli and H.~Gomez and G.~Lorenzo and G.~Marinoschi 
                  and A.~Reali and E.~Rocca},
  TITLE =        {Optimal control of cytotoxic and antiangiogenic
                  therapies on prostate cancer growth},
  JOURNAL =      {Math. Models Methods in Appl. Sci.},
  YEAR =         {2021},
  volume =       {31},
  pages =        {1419--1468},
}

\end{document}